\documentclass{article}
\usepackage[a4paper, total={16.5cm, 26cm}]{geometry}
\usepackage[table]{xcolor}

\usepackage[scr=rsfs]{mathalpha}

\usepackage{tikz}
\usetikzlibrary{}
\usetikzlibrary{matrix,chains,positioning,decorations.pathreplacing,arrows}
\usetikzlibrary{positioning,calc}
\usetikzlibrary{graphs,graphs.standard,quotes}

\usepackage{pgfplots}
\pgfplotsset{compat=1.17} 
\usepackage[utf8]{inputenc}
\usepackage[square,sort,comma,numbers]{natbib}
\usepackage{graphicx}
\graphicspath{ {./images/} }
\usepackage{amsmath}
\usepackage{amsfonts}
\usepackage{amssymb}
\usepackage{algorithm}
\usepackage{algpseudocode}
\usepackage{mathtools}
\usepackage{mathrsfs}
\usepackage{tabularx}
\usepackage{multirow}
\usepackage{float}
\usepackage{subcaption}
\usepackage[font={footnotesize}, margin=0.5cm]{caption}
\usepackage{stackengine}
\usepackage{url}
\usepackage{booktabs}
\usepackage{makecell}
\usepackage{lineno}
\usepackage{bbm}
\usepackage[bookmarks,hidelinks]{hyperref}

\newcommand{\abs}[1]{\left|#1\right|}
\newcommand{\norm}[1]{\left\|#1\right\|}
\newcommand{\normdg}[1]{{\left\vert\kern-0.25ex\left\vert\kern-0.25ex\left\vert #1 
    \right\vert\kern-0.25ex\right\vert\kern-0.25ex\right\vert}}

\newcommand{\bref}[1]{(\ref{#1})}
\newcommand{\jmp}[1]{[\![#1]\!]}
\newcommand{\avg}[1]{\{\!\!\{#1\}\!\!\}}

\algrenewcommand\algorithmiccomment[2][\footnotesize]{{#1\hfill\(\triangleright\) \textit{#2}}}

\newlength{\bibitemsep}
\newlength{\bibparskip}
\let\oldthebibliography\thebibliography
\renewcommand\thebibliography[1]{%
  \footnotesize
  \oldthebibliography{#1}%
  \setlength{\parskip}{\bibitemsep}%
  \setlength{\itemsep}{\bibparskip}%
}

\title{Discovering Artificial Viscosity Models for Discontinuous Galerkin approximation of Conservation Laws using Physics-Informed Machine Learning}
\author{Matteo Caldana$^{a,}$\thanks{Corresponding author: {\tt matteo.caldana@polimi.it}},
Paola F. Antonietti$^a$, 
Luca Dede'$^a$\\[0.3cm]
\small\textit{$^a$MOX, Dipartimento di Matematica, 
Politecnico di Milano,
Piazza Leonardo da Vinci 32,
20133 Milano, Italy
}}

{\date{\small\today}}

\begin{document}

\maketitle
{\subsection*{\centering Abstract}
\small 
\begin{changemargin}{1cm}{1cm}
    Finite element-based high-order solvers of conservation laws offer large accuracy but face challenges near discontinuities due to the Gibbs phenomenon. Artificial viscosity is a popular and effective solution to this problem based on physical insight. In this work, we present a physics-informed machine learning algorithm to automate the discovery of artificial viscosity models. \textcolor{black}{We refer to the proposed approach as an ``\textit{hybrid}'' approach which stands at the edge between supervised and unsupervised learning. More precisely, the proposed ``\textit{hybrid}'' paradigm is not supervised in the classical sense as it does not utilize labeled data in the traditional way but relies on the intrinsic properties of the reference solution}. The algorithm is inspired by reinforcement learning and trains a neural network acting cell-by-cell (the viscosity model) by minimizing a loss defined as the difference with respect to a reference solution thanks to automatic differentiation. This enables a dataset-free training procedure. We prove that the algorithm is effective by integrating it into a state-of-the-art Runge-Kutta discontinuous Galerkin solver. We showcase several numerical tests on scalar and vectorial problems, such as Burgers' and Euler's equations in one and two dimensions. Results demonstrate that the proposed approach trains a model that is able to outperform classical viscosity models. Moreover, we show that the learnt artificial viscosity model is able to generalize across different problems and parameters.
\end{changemargin}
}
\vspace{0.2cm}
\noindent\textbf{Key words:} Artificial viscosity, Conservation laws, Discontinuous Galerkin, Physics-informed machine learning, Neural networks, Reinforcement learning.

\noindent\textbf{MSC subject classification:} 35L65, 65M60, 68T01

\section{Introduction}

Hyperbolic conservation laws, such as Euler equations or the equations of magnetohydrodynamics, have a wide range of applications in Engineering, ranging from aerodynamics, weather predictions, and civil engineering to astrophysics and nuclear fusion. The Discontinuous Galerkin (DG) method \cite{cockburn1989tvb, cockburn2012discontinuous, dumbser2006arbitrary, hesthaven2007nodal, di2011mathematical} is a well-established (and continuously growing in popularity) high-order solver that offers conservation of physical properties, high-order accuracy, geometric flexibility \cite{antonietti2016review}, and parallel efficiency \cite{klockner2009nodal, landmann2008parallel}. A known drawback of any high-order methods is that in the presence of discontinuities of the solution, the Gibbs phenomenon may appear \cite{gottlieb1997gibbs}, which causes spurious numerical oscillations that invalidate the physical meaning of the solution. Indeed, Godunov's theorem \cite{wesseling2009principles} states that any linear monotone numerical scheme can be at most first-order accurate. This problem is especially relevant since in certain regions of the domain continuous initial conditions may lead to a loss of regularity of the solution and thus numerical instability appears.

Several approaches have been developed in the literature to correct and reduce the effect of the Gibbs phenomenon for DG methods. A common feature of most of these approaches is that they rely on the identification of troubled cells, that is mesh elements where the solution loses regularity. For instance, slope limiters \cite{cockburn1998runge, burbeau2001problem, kuzmin2010vertex} are widely used, but they can significantly compromise the global solution accuracy, and an improper parameter selection may lead to increased computational costs. Other authors have investigated Weighted Essentially Non-Oscillatory (WENO) reconstruction \cite{shu2003high, qiu2004hermite, qiu2005runge}, which maintains high-order accuracy but may exhibit a large computational cost. More recently, the Multidimensional Optimal Order Detection (MOOD) approach \cite{dumbser2014posteriori} and filtering methods \cite{orlando2023filtering} have been proposed. These approaches switch the employed method on the troubled cells: from a high-order DG method to a monotone first-order method. However, their performance relies on the quality of the indicator and parameter tuning, respectively.

Instead, in this paper, we focus on a different approach called Artificial Viscosity (AV), which is based on the addition of a suitable amount of dissipation near discontinuities so to recover a smooth solution and thus high-order accuracy. A common characteristic of artificial viscosity approaches is the definition of a metric to evaluate the regularity of the solution. For instance, they may use differential operators that have a particular physical meaning \cite{mani2009suitability}, the jump of the flux on the element boundaries \cite{bassi2009high}, or the residual in the mesh element \cite{hartmann2006adaptive}. Other models instead rely on estimating the average decay rate of modal coefficients \cite{klockner2011viscous}, or just targeting the highest modes \cite{persson2006sub}. A very successful model is the so-called entropy viscosity (EV) model, which employs as a measure of regularity the residual of an entropy pair \cite{guermond2011entropy, zingan2013implementation}. A significant drawback of all these models is that their accuracy, robustness, and stability strongly depend on parameters tuned by physical insight and intuition. Indeed, the lack of an established rule for estimating optimal values and the fact that parameters must be picked on a problem-dependent basis makes their usage expensive.

In recent years, a large effort has been poured into overcoming these problems by employing supervised learning methods and more specifically, deep learning. In particular, neural networks (NN) have been used to create a universal artificial viscosity model that does not require to be tuned on a problem-dependent basis. Indeed, in \cite{discacciati2020controlling, schwander2021controlling} the authors prove that is possible to train a NN that identifies regions where dissipation of numerical oscillations is needed. In \cite{tassi2023machine} the authors use a NN to determine the optimal parameter for the SUPG stabilization method. In general, NNs have been proven to be a valuable tool to surrogate models and automate the expensive phase of trial-and-error present in many numerical applications \cite{heinlein2021combining, eichinger2020stationary, fresca2021comprehensive, caldana2023deep, ray2019detecting}.

Supervised learning however has some limitations. On one hand, building an extensive and high-quality dataset demands a considerable investment of time and computational resources. Indeed, in this case, its construction entails an expensive phase of parameter tuning of classical viscosity model in order to create examples that the NN could mimic. Moreover, supervised learning is not able to generalize beyond what is present in the dataset: the NN may generalize to problems not present in the dataset, but it is not possible to surpass the accuracy of the best among the classical viscosity models, which do not guarantee optimality. Namely, a critical issue of this kind has been exposed in a-posteriori testing of turbulence models trained by supervised learning \cite{wu2018physics}. A possible solution proposed in \cite{novati2021automating} is to employ reinforcement learning (RL) \cite{sutton2018reinforcement}. RL is a well-established but expensive paradigm where an agent learns to make decisions by interacting with an environment, receiving feedback in the form of rewards, and adjusting its behavior to maximize cumulative reward over time. 

Our goal in this paper is to propose and test an algorithm to train a NN-based artificial viscosity model in a \textcolor{black}{\textit{hybrid}} learning framework\textcolor{black}{, which stands at the edge between supervised and unsupervised learning}. To this end, we formulate the problem as an RL task and we enhance it with physics-informed machine learning (PIML) \cite{karniadakis2021physics, hao2022physics}. Namely, thanks to automatic differentiation, we can differentiate through the RL environment, effectively embedding information on the physics of the problem into the learning dynamics. The resulting approach shares some similarities with model-based RL such as \cite{hafner2019dream}, differentiable RL \cite{lutter2021differentiable}, and optimal control. Therefore, our approach overcomes the shortcomings of supervised learning and is cheaper than RL since we can exploit the differentiability of the environment (the DG solver). The reward is defined as the dissimilarity from a reference solution with the addition of suitable regularization terms. This also enables seamless integration of noisy data obtained from measurements of quantities of interest into the loss function, thereby being able to learn new physical models directly from data. The NN that we employ is a generalization of the one presented in \cite{discacciati2020controlling} aimed at increasing flexibility and accuracy. An ablation study proves the efficacy of the proposed modifications.

\vspace{2mm}
The paper is organized as follows. Section~\ref{sec:math-model} introduces the mathematical model of conservation laws and its semi-discrete and algebraic formulations. In Section~\ref{sec:av}, we define the classical artificial viscosity models and we introduce the NN viscosity model. In Section~\ref{sec:physics-informed-machine-learning} we introduce the novel training algorithm we are proposing and discuss the advantages of physics-informed machine learning. In Section~\ref{sec:results} we show convergence results for a smooth problem and we present several applications to problems with discontinuous solutions. Finally, in Section~\ref{sec:conclusions}, we draw some conclusions and discuss further developments.

\section{The mathematical model}\label{sec:math-model}
In this section, we state the formulation of conservation law under consideration. The problem is
defined in a bounded open domain in space $\boldsymbol x \in \Omega \subset \mathbb{R}^d$ $(d = 1, 2)$ and for a time $t \in (0, T]$, with $T>0$ given. The unknown of our problem is the vector of $m$ conserved variables $\boldsymbol u = \boldsymbol u(\boldsymbol x, t) : \Omega \times (0, T] \to \mathbb{R}^m$. The considered convention-diffusion problem reads as follows:
\begin{equation}
    \begin{cases}
    \partial_t \boldsymbol u + \nabla \cdot (- \mu(\boldsymbol u) \nabla \boldsymbol u + \boldsymbol f(\boldsymbol u))  = 0 \quad &\text{in} \: \Omega \times (0, T] ,\\
    \boldsymbol u(\boldsymbol x, 0) = \boldsymbol u_0(\boldsymbol x) &\text{in} \: \Omega \times \{0\}
    \end{cases}
    \label{eq:strong_conservation_law}
\end{equation}
endowed with suitable boundary conditions on $\partial\Omega \times (0, T]$. Namely, we either employ periodic boundary conditions or mixed Dirichlet-Neumann boundary conditions

\begin{equation*}
    \begin{cases}
        \boldsymbol u(\boldsymbol{x}, t) = \boldsymbol{g}_D(\boldsymbol{x}) &\text{on} \: \Gamma_D \times (0,T],\\
        \boldsymbol{f}(\boldsymbol u(\boldsymbol{x}, t)) \boldsymbol{n} = \boldsymbol{g}_N(\boldsymbol{x}) &\text{on} \: \Gamma_N \times (0,T],
    \end{cases}
\end{equation*}
where $\Gamma_N$ and $\Gamma_D \neq \emptyset$ are disjoint sets $\mathring{\Gamma}_N \cap \mathring{\Gamma}_D = \emptyset$ such that $\partial\Omega = \overline{\Gamma}_N \cup \overline{\Gamma}_D$ and $\boldsymbol{g}_D : \Gamma_D \to \mathbb{R}^m, \boldsymbol{g}_N : \Gamma_N \to \mathbb{R}^m$ are a regular enough boundary datum. In Eq.~\bref{eq:strong_conservation_law}, $\boldsymbol f  : \mathbb{R}^m \to \mathbb{R}^{m \times d}$ is the advection flux and $\boldsymbol u_0: \Omega \to \mathbb{R}^m$ is a given initial condition. The viscous flux is controlled by the viscosity coefficient $\mu : \mathbb{R}^m \to \mathbb{R}$ that is a positive and bounded function that has the physical meaning of an artificial diffusion term.

\subsection{The semi-discrete discontinuous Galerkin formulation}
For the spatial discretization, we employ the discontinuous Galerkin method. Let $\mathcal K_h$ be a shape-regular family of open, disjoint cells (intervals for $d=1$ or triangles for $d=2$) $K$ that approximates $\Omega$. Moreover, let $\mathcal{E}_h$ be the set of all faces (points for $d=1$ or edges for $d=2$). On each face $F\in\mathcal{E}_h$, approximating the boundary $\partial\Omega$, a single type of boundary condition is imposed at a time (e.g.\ only Dirichlet or Neumann).
We define $h = \max_{K \in \mathcal K_h} h_K$, where $h_K$ denotes the diameter of $K$. 
We define the averages on and \textcolor{black}{the} jumps across the edges for a smooth enough scalar function $v:\Omega\to\mathbb{R}$\textcolor{black}{,} for a smooth enough vector-valued function $\boldsymbol q:\Omega\to\mathbb{R}^n$ \textcolor{black}{ and a smooth enough tensor-valued function $\boldsymbol \psi$}. Let us consider two neighboring elements $K^+$ and $K^-$, let $\boldsymbol n^+$ and $\boldsymbol n^-$ be their outward normal unit vector and let $v^\pm, \boldsymbol q^\pm\textcolor{black}{, \boldsymbol{\psi}^\pm}$ be the restrictions to $K^\pm$, respectively. In accordance with \cite{arnold2002unified}, \textcolor{black}{defined the symmetric outer product $\boldsymbol q \odot \boldsymbol n = \frac{1}{2} (\boldsymbol v \otimes \boldsymbol n + \boldsymbol n \otimes \boldsymbol v)$,} we define

\begin{itemize}
    \item on the set of interior edges (denoted as $\mathring{\mathcal{E}}_h$):
    \begin{gather*}
        \avg{v} = \frac{1}{2}(v^- + v^+), \qquad 
        \jmp{v} = v^+ \boldsymbol n^+ + v^- \boldsymbol n^-,\\
        \avg{\boldsymbol q} = \frac{1}{2}(\boldsymbol q^- + \boldsymbol q^+), \qquad
        \jmp{\boldsymbol q} = \boldsymbol q^+ \odot \boldsymbol n^+ + \boldsymbol q^- \odot \boldsymbol n^-,\\
        \textcolor{black}{\avg{\boldsymbol{\psi}} = \frac{1}{2}(\boldsymbol \psi^- + \boldsymbol \psi^+),} \qquad
        \textcolor{black}{\jmp{\boldsymbol{\psi}} = \boldsymbol \psi^- \boldsymbol n^- +  \boldsymbol \psi^+ \boldsymbol n^+.}
    \end{gather*}
    \item on boundary edges\textcolor{black}{:}
    \begin{gather*}
        \avg{v} = v, \quad \avg{\boldsymbol q} = \boldsymbol q, \quad \textcolor{black}{\avg{\boldsymbol \psi} = \boldsymbol \psi, \quad}\\
        \jmp{v} = v \boldsymbol n, \quad \jmp{\boldsymbol q} = \boldsymbol q \odot \boldsymbol n, \quad \textcolor{black}{\jmp{\boldsymbol \psi} = \boldsymbol \psi \boldsymbol{n}}.
    \end{gather*}
\end{itemize}

To further stabilize the scheme, as shown in \cite{brezzi2004discontinuous, hesthaven2007nodal}, we \textcolor{black}{introduce} the following weighted average, which introduces a stabilizing numerical flux
\begin{equation}
    \avg{\boldsymbol q}_\alpha = \avg{\boldsymbol f(\boldsymbol q)} + \frac{\textcolor{black}{\alpha}}{2}\jmp{\boldsymbol q} \quad \forall\;F\in\mathcal{E}_h,
    \label{eq:flux}
\end{equation}
where $\alpha \, \textcolor{black}{> 0}$ \textcolor{black}{is} a scalar. For some choices of $\alpha$ we can recover already well-known numerical fluxes: if $\boldsymbol f$ is the identity and $\textcolor{black}{\alpha = 1}$ we have the upwind flux \cite{ayuso2009discontinuous}, if $\alpha$ is the largest eigenvalue (in absolute value) $\lambda$ of the Jacobian $\frac{\partial \boldsymbol f}{\partial \boldsymbol u}$ we have the Rusanov flux \cite{rusanov1962calculation}.
To construct the semi-discrete formulation, we define for any integer $k \geq 1$ the finite element space of discontinuous piecewise polynomial functions as
\begin{equation*}
    X_k^h = \{v \in L^2(\Omega) : v|_K \in \mathbb{P}^k(K) \;\; \forall \; K \in \mathcal{K}_h\}, \quad \boldsymbol V_k^h = [X_k^h]^d,
\end{equation*}
where $\mathbb{P}^k(K)$ is the space of polynomials of total degree less than or equal to $k$ on $K$.
Hence, the semi-discrete formulation reads as follows:
For any $t \in (0, T ]$, find $\boldsymbol u_h(t) \in \boldsymbol V_k^h$ such that:
\begin{equation}
    \begin{cases}
        (\dot{\boldsymbol{u}}_h(t), \boldsymbol v_h)_\Omega + \mathscr{A}_h(\boldsymbol u_h(t), \boldsymbol v_h) - \mathscr{B}_h(\boldsymbol u_h(t), \boldsymbol v_h) + \mathscr{I}_h(\boldsymbol u_h(t), \boldsymbol v_h) = 0 \qquad &\forall \; \boldsymbol v_h \in \boldsymbol V_k^h,\\
        \boldsymbol u_h(0) = \boldsymbol u_{0h} &\text{in} \: \Omega \times \{0\}, 
    \end{cases}
    \label{eq:weak_form}
\end{equation}
where
\begin{itemize}
    \item $\mathscr{A}_h: \boldsymbol V_k^h \times \boldsymbol V_k^h \rightarrow \mathbb R$ is defined by
    \begin{equation*}
        \begin{aligned}
            \mathscr{A}_h(\boldsymbol u_h, \boldsymbol v_h) &= \sum_{K\in\mathcal{K}_h}\int_\Omega \mu(\boldsymbol u_h) \nabla \boldsymbol u_h \textcolor{black}{:} \nabla \boldsymbol v_h \textcolor{black}{\, \mathrm d \Omega} \\
            &+ \sum_{F \notin \Gamma_N}\frac{\tau k^2}{\abs{F}}\int_{F}\mu(\boldsymbol u_h)\jmp{\boldsymbol u_h} \textcolor{black}{:} \jmp{\boldsymbol v_h} \textcolor{black}{\, \mathrm d \Sigma}\\
            &-\sum_{F \in \mathring{\mathcal{E}}_h} \textcolor{black}{\int_{F}} \avg{\nabla \boldsymbol u_h} \textcolor{black}{:} \jmp{\boldsymbol v_h}+\jmp{\boldsymbol u_h} \textcolor{black}{:} \avg{\nabla \boldsymbol v_h} \textcolor{black}{\, \mathrm d \Sigma} \quad \forall \; \boldsymbol u_h, \boldsymbol v_h \in \boldsymbol V_k^h,\\ 
        \end{aligned}
    \end{equation*}
    where $\tau > 0$ is a suitable penalization constant and $\abs{F}$ is the measure in the topology of the boundary.
    \item $\mathscr{B}_h: \boldsymbol V_k^h \times \boldsymbol V_k^h \rightarrow \mathbb R$ is defined as
    \begin{equation*}
      \begin{aligned}
        \mathscr{B}_h(\boldsymbol u_h, \boldsymbol v_h) &= \sum_{K \in \mathcal{K}_h} \int_{K} \boldsymbol{f}(\boldsymbol u_h) \textcolor{black}{:} \nabla \boldsymbol v_h \textcolor{black}{\, \mathrm d \Omega}\\ 
        & - \sum_{F \in \mathring{\mathcal{E}}_h}\int_{F} \avg{\boldsymbol u_h}_\lambda \textcolor{black}{:} \jmp{\boldsymbol v_h} \textcolor{black}{\, \mathrm d \Sigma} \quad \forall \; \boldsymbol u_h, \boldsymbol v_h \in \boldsymbol V_k^h,
      \end{aligned}
    \end{equation*}
    \textcolor{black}{where $\avg{\boldsymbol u_h}_\lambda$ indicates the Rusanov flux, as defined in Eq.~\bref{eq:flux}.}
    \item $\mathscr{I}_h: \boldsymbol V_k^h \times \boldsymbol V_k^h \rightarrow \mathbb R$ is a form that includes suitable terms for the boundary conditions.
\end{itemize}

\subsection{Algebraic formulation}\label{sec:algebraic-form}
Denote as $(\boldsymbol \varphi_i)_{i=1}^I$ a suitable basis for the discrete space $\boldsymbol V_k^h$, where $I=\text{dim}(\boldsymbol V_k^h)$, then
\begin{equation*}
    \boldsymbol u_h(t) = \sum_{i=1}^I U_i(t) \boldsymbol \varphi_i.
\end{equation*}
We collect the coefficients $U_i(t)$ of the expansions of $\boldsymbol{u}_h(t)$ in the vector function $\boldsymbol U: (0, T] \to \mathbb R^{mI}$. By defining the mass matrix
\begin{equation*}
    [M]_{ij} = (\boldsymbol{\varphi}_i, \boldsymbol{\varphi}_j)_\Omega, \quad \forall\; i,j=1,...,I
\end{equation*}
and the non-linear terms
\begin{equation*}
    [A(\boldsymbol{U}(t))]_{i} = \mathscr{A}_h(\boldsymbol{u}_h, \boldsymbol{\varphi}_i), \quad 
    [B(\boldsymbol{U}(t))]_{i} = \mathscr{B}_h(\boldsymbol{u}_h, \boldsymbol{\varphi}_i), \quad 
    [I(\boldsymbol{U}(t))]_{i} = \mathscr{I}_h(\boldsymbol{u}_h, \boldsymbol{\varphi}_i), \quad \forall\; i=1,...,I.
\end{equation*}
We rewrite problem \bref{eq:weak_form} in algebraic form as follows

\begin{equation}
    \begin{cases}
        M\dot{\boldsymbol{U}}(t) + A(\boldsymbol{U}(t)) - B (\boldsymbol{U}(t)) + I(\boldsymbol{U}(t))= 0 \qquad t \in (0, T],\\
        \boldsymbol{U}(0) = \boldsymbol U_0,
    \end{cases}
    \label{eq:algebraic_form}
\end{equation}
where $\boldsymbol U_0$ is the projection of the initial datum $\boldsymbol{u}_0$ on $\boldsymbol V_k^h$. For the space discretization, we employ a nodal spectral basis derived from a multidimensional generalization of the Legendre-Gauss-Lobatto points as shown in \cite{hesthaven2007nodal}.

\subsection{Fully-discrete formulation}
System \bref{eq:algebraic_form} can be solved with various time-marching strategies, like Runge-Kutta schemes, by defining a partition of $N$ intervals $0 = t^0 <
t^1 < ... < t^N = T$. In literature (see \cite{hesthaven2007nodal}), this is usually coupled with an adaptive time-step $\Delta t = \Delta t(t)$ determined by the following CFL condition at a time $t$:
\begin{equation}
    \Delta t = \dfrac{\textnormal{CFL}}{\displaystyle\frac{k^2}{h} \max_{\boldsymbol{x}\in\Omega}\abs{\boldsymbol f'(\boldsymbol u_h(\boldsymbol{x}, t))}+\frac{k^4}{h^2} \max_{\boldsymbol{x}\in\Omega}\mu(\boldsymbol u_h(\boldsymbol{x}, t))}.
    \label{eq:cfl}
\end{equation}
Specifically, we utilize one of the two following options: a third-order Strong Stability Preserving Runge-Kutta (SSPR3) scheme \cite{durran2010numerical} or a fourth-order low-storage explicit Runge-Kutta with five stages \cite{cockburn1989tvb}. The choice between the two is made so to minimize the number of total stages while maintaining a stable solution. \textcolor{black}{To further verify the consistency of the results, we compute the maximum Courant number}
\begin{equation*}
    \textcolor{black}{C = \max_{\boldsymbol{x}\in\Omega}\abs{\boldsymbol f'(\boldsymbol u_h(\boldsymbol{x}, t))} \frac{\Delta t}{h}.}
\end{equation*}

\subsection{Convective flux models}
This study investigates various conservation laws. For one-dimensional ($d=1$) scenarios, we examine
\begin{itemize}
    \item the linear advection equation ($m=1$): 
    \begin{equation}
        f(u) = \beta u,
        \label{eq:adv-1d}
    \end{equation}
    where the transport field $\beta \in \mathbb{R}$ is constant;
    \item the Burgers' equation ($m=1$):
    \begin{equation}
        f(u) = \frac{1}{2}u^2;
        \label{eq:burgers-1d}
    \end{equation}
    \item the Euler system, namely, given the ideal gas law for pressure $p = (\gamma - 1)(e - \frac{1}{2}\rho\abs{v}^2)$, where $\gamma = 7 / 5$ is the heat capacity ratio, $e$ is the \textcolor{black}{total} energy, $\rho$ is the density and $v$ is the velocity, we have ($m=3$):
    \begin{equation}
        \boldsymbol u = 
            \begin{bmatrix}
            \rho\\
            \rho v\\
            e
            \end{bmatrix},\qquad
            \boldsymbol f(\boldsymbol u) =
            \begin{bmatrix}
                \rho v\\
                \rho v^2 + p\\
                 v(e + p)
            \end{bmatrix}.
        \label{eq:euler-1d}
    \end{equation}
\end{itemize}
The two-dimensional ($d=2$) models we consider are
\begin{itemize}
    \item the linear advection equation ($m=1$): 
    \begin{equation}
        \boldsymbol f(u) = \boldsymbol{\beta} u,
        \label{eq:advection-2d}
    \end{equation}
    where the transport field $\boldsymbol{\beta} \in \mathbb{R}^2$ is constant;
    \item the KPP rotating wave problem \cite{guermond2011entropy,kurganov2007adaptive} ($m = 1$), characterized by a non-convex flux in each component 
    \begin{equation}
        \boldsymbol f(u) = \begin{bmatrix}\sin u\\ \cos u\end{bmatrix};
        \label{eq:KPP-2d}
    \end{equation}
    \item the Euler system ($m=4$), defined as 
    \begin{equation}
        \boldsymbol u = 
            \begin{bmatrix}
            \rho\\
            \rho v_1\\
            \rho v_2\\
            e
            \end{bmatrix},\qquad
            \boldsymbol f(\boldsymbol u) =
            \begin{bmatrix}
                \rho v_1 & \rho v_2\\
                \rho v_1^2 + p & \rho v_2^2 + p\\
                \rho v_1v_2 & \rho v_1v_2\\
                 v_1(e + p) & v_2(e + p)
            \end{bmatrix},
            \label{eq:euler-2d}
    \end{equation}
    where $\boldsymbol v = \begin{bmatrix}v_1\\ v_2\end{bmatrix}$ is the velocity field, \textcolor{black}{$e$ is the total energy}, and the pressure $p = (\gamma - 1)(e - \frac{1}{2}\rho\norm{\boldsymbol{v}}_2^2)$ is defined analogously as before.
\end{itemize}

\section{Artificial viscosity}\label{sec:av}
We provide a concise review of classical artificial viscosity models that hold a pivotal role in understanding the limitations of current methods. We then present the EV model which is the state-of-the-art model serving as the main tool for validating the proposed deep learning physics-informed technique. Finally, we introduce the artificial viscosity models based on neural networks. For the sake of simplicity, we describe the procedure for scalar problems. The standard way to generalize to systems is to choose one of the variables as the representative scalar quantity to be used within the viscosity model.

\subsection{An outline of artificial viscosity}
\label{sec:av-outline}
Most artificial viscosity models are based on a heuristic that may be very different from model to model. 
However, it is still possible to recognize some common features that all artificial viscosity models share:
\begin{itemize}
    \item For each element $K\in\mathcal{K}_h$, a shock indicator measures the size of the local discontinuities and outputs a viscosity $\mu_K$ to be added in that cell. The viscosity $\mu_K$ must be small or null when the solution is smooth, to avoid modifications to the physics of the problem. In the presence of discontinuities, $\mu_K$ must be the smallest value that is sufficient to eliminate Gibbs oscillations. 
    \item Since large viscosity values produce over-dissipation, rendering the model unable to capture significant features of the solution, at every time instant $t$, and for each element $K \in \mathcal{K}_h$ we define a viscosity threshold
    \begin{equation}
        \mu_K^\text{bnd}(\boldsymbol{x}, t) = \max \{\max_{\boldsymbol{x} \in K } \mu_K(\boldsymbol{x}, t), \mu_{\max{}}(t)\}, \quad \mu_{\max{}}(t) = c_{\max{}} \frac{h}{k} \max_{\boldsymbol{x} \in K }{\abs{\boldsymbol{f}'(u_h(\boldsymbol{x}, t))}},
        \label{eq:max_visc}
    \end{equation}
    where $c_{\max{}}$ is a scalar parameter to be tuned for each problem \cite{persson2006sub,klockner2011viscous} and $\boldsymbol f'$ defines the velocity of propagation of the problem. For $m>1$ we consider the restriction to the equation of the representative variable. \textcolor{black}{This restriction does not constitute a limitation; all the components of $\boldsymbol f'$ can be utilized as inputs for the neural network. The primary drawback of this approach is the increased computational cost.} Notice that the definition of $\mu_{\max{}}$ assures that the viscous contribution is at most linear.
    \item The viscosity is globally smoothed as outlined in \cite{yu2020study}. Despite the existence of advanced techniques \cite{abbassi2014shock}, we focus on a straightforward three-step piecewise linear interpolator approach (this procedure can easily be generalized to an arbitrary degree):
    \begin{enumerate}
        \item On each vertex of the discretization $\mathcal K_h$, compute the average viscosity among all the elements that share that vertex;
        \item Compute the element-wise linear interpolator, given the average viscosities on the vertexes of each element;
        \item Evaluate the interpolator on the required points in order to build the algebraic system.
    \end{enumerate}
    Our numerical findings and \cite{barter2010shock} support the thesis that the smoothing mitigates numerical oscillations. 
\end{itemize}

In this study, we consider updates to the viscosity which are computed at every time step. However, it is possible to update the viscosity less frequently in order to decrease the computational cost.

\subsection{Entropy viscosity}
The entropy viscosity model \cite{guermond2011entropy, zingan2013implementation, dafermos2005hyperbolic} is a state-of-the-art model that produces a viscosity coefficient based on the local size of entropy production. 
Namely, it is known \cite{bardos1979first} that the scalar inviscid continuous initial boundary value problem \bref{eq:strong_conservation_law} has an unique entropy solution satisfying 
\begin{equation*}
    \partial_t E(u) + \nabla \cdot \boldsymbol F(u) \leq 0
\end{equation*}
where the two functions $E:\mathbb{R} \to \mathbb{R}$ and $F:\mathbb{R}\to\mathbb{R}^d$ are such that
\begin{equation*}
    \boldsymbol F(w) = \int E'(w)\boldsymbol f'(w) \;\textnormal{d}w, \qquad E \textnormal{ is convex, }
\end{equation*}
and are called an entropy pair.
We then define the local cell viscosity as
\begin{equation}
    \mu_K(\boldsymbol{x}, t) = c_K \left( \frac{h}{k} \right)^2 \frac{\max\{\abs{D_h(\boldsymbol{x}, t)}, \abs{H_h(\boldsymbol{x}, t)}\}}{\norm{E(u_h) - \overline{E}(u_h)}_\infty}
    \label{eq:ev-visc}
\end{equation}
where $c_K$ is a parameter that must be manually tuned for each problem, $\overline{E}$ is the spatial integral average of the entropy over the domain, $D_h$ is the residual of the space-time entropy equation
\begin{equation*}
    D_h(\boldsymbol{x}, t) = \partial_t E(u_h(\boldsymbol{x}, t)) + \nabla \cdot \boldsymbol F(u_h(\boldsymbol{x}, t)) ,
\end{equation*}
and $H$ is the effect of the jump of the entropy flux across the element boundaries
\begin{equation*}
    H_h(\boldsymbol{x}, t) = \left(\frac{h}{k}\right)^{-1} \jmp{\boldsymbol F(u_h(\boldsymbol{x}, t))} \cdot \boldsymbol n.
\end{equation*}
Finally, we impose an upper limit as prescribed by \bref{eq:max_visc}. Let us now use this model to make some observations that actually extend to almost all AV models. In particular, we can outline three main weaknesses of classical AV models that we would like to overcome:
\begin{enumerate}
    \item [(A)] The evaluation of the viscosity $\mu_K^\text{bnd}$ may be computationally expensive. In this particular case, computing $D_h$ and $H_h$ is expensive.
    \item [(B)] They depend on parameter(s) that must be chosen manually. Specifically, $c_{\max{}}$ and $c_K$ have to be tuned with trial-and-error. This requires deep know-how of the model or a time-consuming phase of trial-and-error.
    \item [(C)] They do not guarantee optimality: e.g., the EV model is based on an entropy heuristic argument and does not provide any error bound estimate.
\end{enumerate}

\subsection{Basic concepts of Neural networks}
Let $N_I, N_O \in \mathbb{N}$ be strictly positive. An artificial neural network is a function $\mathcal{F} \, :\, \mathbb{R}^{N_{I}} \rightarrow \mathbb{R}^{N_{O}}$ that maps an input vector $\boldsymbol{x} \in \mathbb{R}^{N_{I}}$ to an output vector $\boldsymbol{y} \in \mathbb{R}^{N_{O}}$ and depends on a set of parameters that can be thought to be ordered into a vector ${\boldsymbol{\theta}}$. In this work, we only consider the simplest version of neural networks: dense feed-forward neural networks (FNN), which are defined as follows. Fix two positive integers $L$ and $Z$. Let $\boldsymbol y^{(0)} = \boldsymbol x$ and $\boldsymbol y^{(L)} = \boldsymbol{y}$. An FNN of depth $L$ and width $Z$ is the composition of $L$ functions called layers defined by
\begin{equation*}
  \boldsymbol y^{(l)} = \sigma^{(l)}(\mathrm W^{(l)} \boldsymbol y^{(l-1)} + \boldsymbol b^{(l)}), \quad l=1,...,L,
\end{equation*}
where $\mathrm{W}^{(l)} \in \mathbb{R}^{N_l \times N_{l-1}}$ are matrices of parameters called weights and $\boldsymbol{b}^{(l)} \in \mathbb{R}^{N_l}$ are vectors of parameters called biases. Here, $N_0=N_I$, $N_L=N_O$ and $N_l = Z$ for $l = 1, ..., L - 1$. Finally, at the end of each layer, a scalar non-linear activation function $\sigma^{(l)}$ is applied component-wise. We employ the same activation function (to be specified later, see Section~\ref{sec:training-in-practice}) for all the layers apart from the last one where we use a softplus activation function to guarantee the positiveness of the output. It can be seen as a smooth approximation of the ReLU, namely:
\begin{equation*}
    \sigma^{(L)}(t) = \log(1 + e^{t})
\end{equation*}

Several learning strategies are available to determine the vector of parameters $\boldsymbol{\theta}$. Among these, one of the most popular paradigms is supervised learning. Namely, we define and subsequently minimize a non-linear cost function, called loss function $\mathcal L$. The optimization is usually performed using a gradient-based optimizer, where the gradient $\nabla_{\boldsymbol{\theta}} \mathcal L$ is computed by means of automatic differentiation \cite{paszke2017automatic}. 
The loss is defined by means of a given dataset of input target pairs $\{(\boldsymbol x, \hat{\boldsymbol y})\}_{i=1}^{N_\textnormal{train}}$ that the NN must learn. Usually, the loss function takes a form similar to the following
\begin{equation}
    \mathcal L(\boldsymbol{\theta}) = \frac{1}{N_\textnormal{train}} \sum_{i = 1}^{N_\textnormal{train}} \norm{\hat{\boldsymbol y} - \boldsymbol{y}(\boldsymbol x;\boldsymbol{\theta})}_{\textcolor{black}{q}},
    \label{eq:loss_classic}
\end{equation}
where $\norm{\cdot}_{\textcolor{black}{q}}$ indicates the discrete $\textcolor{black}{q}$-norm, $\textcolor{black}{q}=1,2$. A sequence of gradient updates made by iterating over the dataset one time is called an epoch. Even if we will propose an alternative paradigm, this is what is used in literature to train the artificial viscosity model and it is useful to fix this as a benchmark to be compared to our approach. 
Regardless of the paradigm used to train the NN, if properly optimized, neural networks are able to generalize, meaning that they make good predictions even for input data not used during the training phase. 
For a comprehensive description of neural networks, we direct the reader to \cite{goodfellow2016deep, zhang2023dive}. 

\subsection{The neural viscosity approach}\label{sec:neural-viscosity}
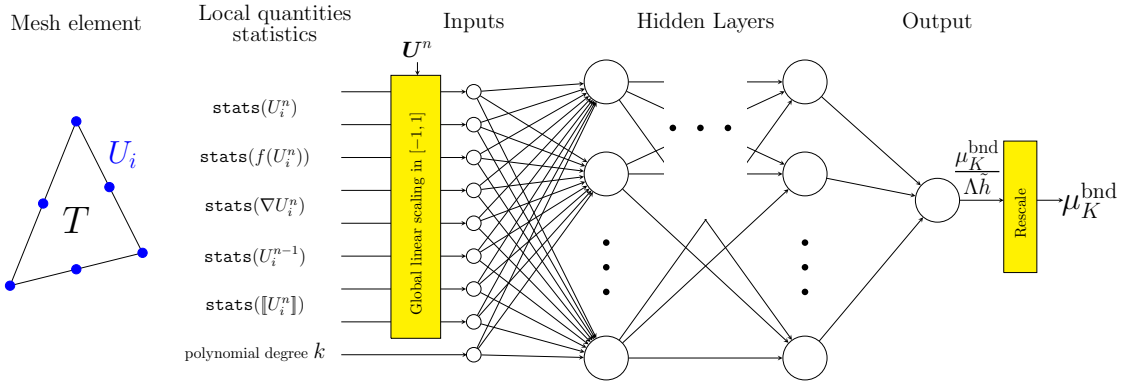
\begin{figure}[t]
    \centering  
    \resizebox{0.9\textwidth}{!}{\tikzset{%
  every neuron/.style={
    circle,
    draw,
    minimum size=1cm
  },
  neuron missing/.style={
    draw=none, 
    scale=4,
    text height=0.333cm,
    execute at begin node=\color{black}$\vdots$
  },
  every input_neuron/.style={
    circle,
    draw,
    minimum size=0.3cm
  },
  input_neuron missing/.style={
    draw=none, 
    scale=4,
    text height=0.333cm,
    execute at begin node=\color{black}$\vdots$
  },
}

\begin{tikzpicture}[x=1.5cm, y=1.5cm, >=stealth]

\foreach \m/\l [count=\y] in {1,2,3,4,5,6,7,8,9}
  \node [every input_neuron/.try, input_neuron \m/.try] (input-\m) at (0,1.95-\y/2) {};

\foreach \m [count=\y] in {1,2,missing,3}
  \node [every neuron/.try, neuron \m/.try ] (hidden1-\m) at (2,3-\y*1.4) {};

\foreach \m [count=\y] in {1,2,missing,3}
  \node [every neuron/.try, neuron \m/.try ] (hidden2-\m) at (5,3-\y*1.4) {};

\foreach \m [count=\y] in {1}
  \node [every neuron/.try, neuron \m/.try ] (output-\m) at (7,0.8-\y) {};

\draw [<-] (input-1) -- ++(-2,0) node [above, midway] {\ };
\draw [<-] (input-2) -- ++(-2,0) node [above, midway] {\ };
\draw [<-] (input-3) -- ++(-2,0) node [above, midway] { \ };
\draw [<-] (input-4) -- ++(-2,0) node [above, midway] { \ };
\draw [<-] (input-5) -- ++(-2,0) node [above, midway] { \ };
\draw [<-] (input-6) -- ++(-2,0) node [above, midway] { \ };
\draw [<-] (input-7) -- ++(-2,0) node [above, midway] { \ };
\draw [<-] (input-8) -- ++(-2,0) node [above, midway] { \ };

\foreach \l [count=\i] in {1,2,3} {
    \ifnum \l=3
        \node [above] at (hidden1-3.north) {};
    \else
        \node [above] at (hidden1-\i.north) {};
    \fi
}

\foreach \l [count=\i] in {1,2,3} {
    \ifnum \l=3
        \node [above] at (hidden2-3.north) {};
    \else
        \node [above] at (hidden2-\i.north) {};
    \fi
}

\foreach \l [count=\i] in {1}
  \draw [->] (output-\i) -- ++(1,0)
    node [above, midway] { \ };

\foreach \i in {1,...,9}
  \foreach \j in {1,...,3}
    \draw [->] (input-\i) -- (hidden1-\j);

\foreach \i in {1,...,3}
  \foreach \j in {1,...,3}
    \draw [->] (hidden1-\i) -- (hidden2-\j);

    \foreach \i in {1,...,3}
    \foreach \j in {1}
    \draw [->] (hidden2-\i) -- (output-\j);

\node[fill=black!0,scale=4,inner xsep=0pt,inner ysep=5mm] at ($(hidden1-1)!.5!(hidden2-2)$) {$\dots$};

\node at (0,2.5)   {\Large Inputs};
\node at (3.5,2.5) {\Large Hidden Layers};
\node at (7,2.5)   {\Large Output};
\node at (-6,2.5)  {\Large Mesh element};
\node[align=center] at (-3,2.5)  {\Large Local quantities\\ \Large statistics};

\node at (-3.3, 1.20) {\large \tt{stats}$(U_i^n)$       };
\node at (-3.3, 0.45) {\large \tt{stats}$(f(U^n_i))$ };
\node at (-3.3,-0.30) {\large \tt{stats}$(\nabla U_i^n)$};
\node at (-3.3,-1.05) {\large \tt{stats}$(U_i^{n - 1})$ }; 
\node at (-3.3,-1.80) {\large \tt{stats}$(\jmp{U^n_i})$ };

\node at (-3.3, -2.5) {polynomial degree \Large $k$};
\draw [<-] (-0.1, -2.55) -- (-2.0, -2.55) node [above, midway] { \ };

\draw[fill=yellow] (-1.25,-2.3) -- (-0.5, -2.3) -- (-0.5, 1.7) -- (-1.25,1.7) -- (-1.25,-2.3);
\node[rotate=90] at (-0.85,-0.4) {Global linear scaling in $[-1, 1]$};
\node at (-0.85, 2.1) {\Large $\boldsymbol{U}^n$};
\draw [<-] (-0.85, 1.7) -- (-0.85, 1.9) node [above, midway] { \ };

\draw[fill=yellow] (8,-1.3) -- (8.5, -1.3) -- (8.5, 0.7) -- (8,0.7) -- (8,-1.3);
\node[rotate=90] at (8.25,-0.4) {Rescale};
\node at (7.6,0.3) {\huge $\frac{\mu_K^\text{bnd}}{\Lambda \tilde{h}}$};
\node at (9.3,-0.2) {\huge $\mu_K^\text{bnd}$};
\draw [->] (8.5,-0.2) -- (8.9,-0.2) node [above, midway] { \ };

\node at (-6,-0.5) {\Huge $T$};
\node[color=blue] at (-5.3,0.5) {\huge $U_i$};

\draw (-7,-1.5) -- (-6, 1) -- (-5, -1) -- (-7,-1.5) ;

\node[mark size=3pt,color=blue] at (-7,-1.5) {\pgfuseplotmark{*}};
\node[mark size=3pt,color=blue] at (-6, 1) {\pgfuseplotmark{*}};
\node[mark size=3pt,color=blue] at (-5, -1) {\pgfuseplotmark{*}};
\node[mark size=3pt,color=blue] at (-6.5, -0.25) {\pgfuseplotmark{*}};
\node[mark size=3pt,color=blue] at (-5.5, 0) {\pgfuseplotmark{*}};
\node[mark size=3pt,color=blue] at (-6, -1.25) {\pgfuseplotmark{*}};

\end{tikzpicture}}
    \caption{Scheme of the online phase of the neural viscosity approach. The neural network returns as output a viscosity value $\mu_K^\text{bnd}$ for each element of the mesh independently of the other elements by extracting relevant statistics \texttt{stats} of the nodal values, namely the minimum, the maximum, the mean, and the standard deviation. Here, $U_i^n$ is the evaluation of $u_h$ at time $t^n$ in the $i$-th degree of freedom of $K$.}
    \label{fig:fnn}
\end{figure}  

The neural viscosity approach, that is the idea of using NNs as a surrogate for AV models, was first proposed in \cite{discacciati2020controlling,schwander2021controlling}. Indeed, NNs are an attractive substitute for classical AV models since they overcome the drawbacks (A) and (B) in the following way:
\begin{enumerate}
    \item [(sA)] NNs are computationally efficient during the online stage since they involve mainly simple matrix-matrix multiplications. Moreover, in recent years large effort was put into developing efficient libraries and hardware to accelerate these linear algebra operations.
    \item [(sB)] NNs can learn highly complex and non-linear functions without the necessity of manually specifying parameters after the training. By leveraging the information present in a large dataset during the training phase, NNs can learn and generalize beyond a single problem.
\end{enumerate}
The architecture that we employ for our NN introduces some improvements with respect to the classical one. The real novelty of this paper lies in how we tackle the third drawback (C), that is how to train the NN and what it learns (or better discovers), which is a matter of Section~\ref{sec:physics-informed-machine-learning}. 

The way neural viscosity works is macroscopically the same in which any AV model works (as described in Section~\ref{sec:av-outline}). The NN works independently on each cell, it receives as inputs some local quantities, and it outputs an artificial viscosity value $\mu_K^\text{bnd}$ for each cell $K$ and finally, it performs a global smoothing step. Let us now describe this procedure in detail. The first difference with respect to the original neural viscosity model is that, instead of building and training a different NN for each polynomial degree $k>0$, we build just one independently of $k$, reducing training costs. This is enabled by the fact that instead of feeding the local nodal values of the solution $\boldsymbol U^n|_K$ (the evaluation of $u_h$ at time $t^n$ in the degrees of freedom of $K$, which is equal to $u_h|_K(t^n)$ since we are employing a nodal basis) to the NN, we extract some representative statistics of these values, namely the mean and the standard deviation and we add as input also the local polynomial degree $k$. The input can be enriched with other information such as the median, the quartiles, the minimum, or the maximum. In Section~\ref{sec:training-in-practice} we present numerical experiments that show the tradeoff between accuracy and computational cost of this choice. The second difference is that to accelerate the training, \textcolor{black}{we employ as inputs of the neural network} not only the current \textcolor{black}{state variable} but \textcolor{black}{all the following features:}
\begin{itemize}
    \item \textcolor{black}{state variable;}
    \item the gradient of the \textcolor{black}{state variable};
    \item the jump of the \textcolor{black}{state variable} across the element boundary, since it is a key information to assess the smoothness of the current solution;
    \item the \textcolor{black}{state variable} at the previous time step, to provide information about temporal correlation;
    \item flux of the \textcolor{black}{state variable}, to provide information about the physics of the problem.
\end{itemize}
Then, as in the original neural viscosity model, we apply a linear scaling to map the inputs in the interval $[-1, 1]$ and a linear rescaling of the targets $\mu_K^\text{bnd}$ of the network. This is common practice for several reasons, including faster convergence during the training (prevents issues like vanishing or exploding gradients) and robustness to variations in input magnitudes (which also is improved generalization). The third difference is that our scaling, instead of being with respect to the local maximum over the element is done with respect to the global maximum. Finally, after the output of the NN, we apply the following rescaling for the element $K$
\begin{equation*}
    y = \frac{\mu_K^\text{bnd}}{\Lambda \tilde{h}}, \quad \textnormal{ where } \quad \Lambda = \max_{K}{\abs{\boldsymbol{f}'(u_h)}}, \quad \tilde{h} = \min\{\max_{\partial K} \abs{\jmp{u_h}}, h\}.
\end{equation*}
Observe that we denote with $y$ the output of the NN since it is a scalar quantity. Moreover, the scaling is time dependent, so it changes at every time step.
The insight behind this rescaling is that $\Lambda$ makes the same network work for \textcolor{black}{different} PDEs since it represents the maximum local wave speed, and $\tilde{h}$ provides the optimal viscosity scaling $h^{k+1}$ in the presence of smooth solution \cite{hesthaven2007nodal}. Indeed, even if the NN prescribes low dissipation when the solution is regular, due to the nature of the softplus function, it is never zero, which prevents the method from having optimal convergence. Instead, $\tilde{h}$ is a parameter-free option that provides the following scaling
\begin{equation*}
    \tilde{h}\sim
    \begin{cases}
        h^{k+1} &\textnormal{ if $u$ is smooth,}\\
        h &\textnormal{ otherwise}
    \end{cases}
\end{equation*} 
which enables a correct reduction of the artificial viscosity when the solution is smooth. \textcolor{black}{In other words, at each timestep, the neural network prescribes a single viscosity value for each element of the mesh.} A complete scheme of what we have just described is shown in Figure~\ref{fig:fnn}.
 
\section{Physics-informed machine learning}\label{sec:physics-informed-machine-learning}

In this section, we describe in detail how the neural network is trained. We propose a novel approach that is inspired by reinforcement learning and is enhanced by physics-informed machine learning. Namely, thanks to automatic differentiation, we can differentiate through the DG solver (the \textit{environment} in RL terminology) and define a loss function (the opposite of the \textit{reward}) as the norm of the error between a reference solution and the solution obtained with the current parameters of the NN (the \textit{policy}) that surrogates a viscosity model (\textit{agent}). The loss is also enhanced with suitable regularization terms. To conclude the comparison with RL terminology, the local statistics extracted from each mesh element are the \textit{state} and the viscosity $\mu_K^\text{bnd}$ the \textit{action} the agent chooses. The major difference with respect to RL is that we are able to differentiate through the environment. Thus, by using a gradient-based optimizer, we can directly compute updates to the parameter of the NN so to minimize the loss (maximize the reward). Hence, this training algorithm enables us to discover a new neural viscosity model that has the property of minimizing the norm of the error. This approach allows us to avoid building a dataset with a reference value of artificial viscosity. Indeed, many RL algorithms overcome the problem of a \textcolor{black}{non-differentiable} environment using a second neural network called \textit{critic} that surrogates it. This is however expensive, since the environment dynamic must also be learnt, and may lead to instabilities due to the interaction of the policy with the critic \cite{lillicrap2015continuous}. \textcolor{black}{The choice of using physics-informed machine learning, and in particular automatic differentiation through the DG solver, is motivated mainly by efficiency. Indeed, automatic differentiation enables the integration of the DG solver dynamics into neural network training, ensuring precision and efficiency. Concerning the architecture, the specific use of a dense feed-forward neural network is motivated by its simplicity, flexibility, and universality. They are straightforward to implement, computationally efficient, and can be tailored to various problems by adjusting their architecture.}

\subsection{Loss function}\label{sec:loss}
\begin{figure}[t]
    \centering  
    \resizebox{0.9\textwidth}{!}{\input{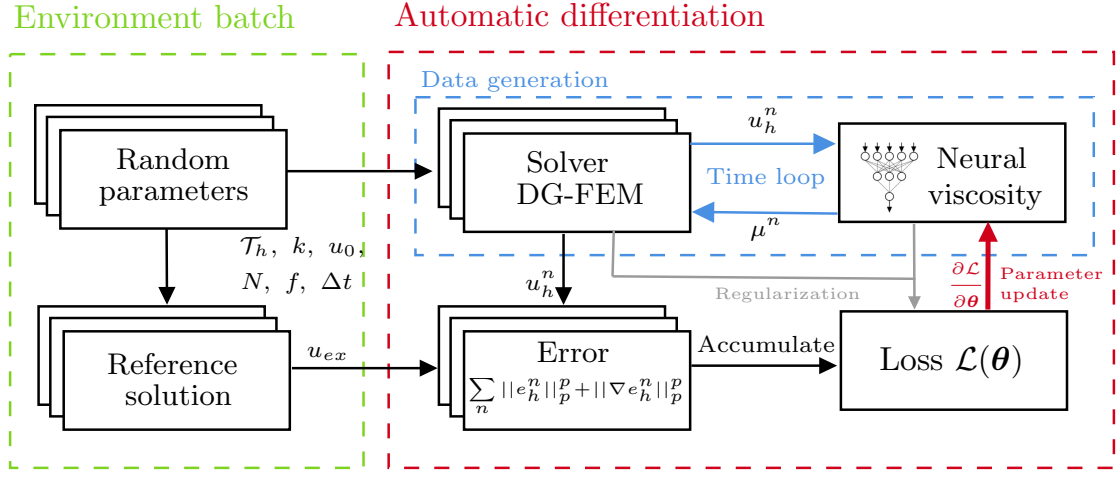}}
    \caption{Schematic description of the proposed training procedure.}
    \label{fig:train}
\end{figure} 

The definition of the loss function is more important in this context than usual since it steers the optimizer toward the discovery of a neural viscosity model with certain properties.
To explain how the loss works, let us make the following remark. Given the discretized initial condition $\boldsymbol U^0$, the solution at the next time step  $\boldsymbol U^1$ depends on the artificial viscosity model we choose, which in return depends on $\boldsymbol U^0$. More in general we could say that 
\begin{equation*}
    \boldsymbol U^n = \boldsymbol U^n(\mu(\boldsymbol U^{n-1}; \boldsymbol{\theta})) = \boldsymbol U^n(\mu(\mu(\boldsymbol U^{n-2}; \boldsymbol{\theta}); \boldsymbol{\theta})) = ... = \boldsymbol U^n(\boldsymbol U^0; \boldsymbol{\theta}), \quad \forall n > 0,
\end{equation*}
where the variable $\boldsymbol{\theta}$ explicits the dependence on the parameters of the NN. Hence, given a discrete reference solution $\boldsymbol U_{\textcolor{black}{ref}}^n$, that is the evaluation of a reference solution $\boldsymbol u_{\textcolor{black}{ref}} : \Omega \times (0, T] \to \mathbb{R}^m$ at the degrees of freedom at time $t^n$, by employ automatic differentiation we can compute
\begin{equation*}
    \frac{\partial}{\partial\boldsymbol{\theta} } \norm{\boldsymbol U^n(\boldsymbol U^0; \boldsymbol{\theta}) - \boldsymbol U_{\textcolor{black}{ref}}^n}_{\textcolor{black}{q}} \quad \forall \: n > 0.
\end{equation*}
This is one of the main ingredients of our loss: using an optimizer, we can find the set of parameter $\boldsymbol{\theta}$ that defines the neural viscosity model that minimizes the error with respect to the reference solution. \textcolor{black}{More precisely, the solution $\boldsymbol{u}_{\textcolor{black}{ref}}$ is either the analytical solution (when available) or an overkill solution obtained using the DG solver on a triangulation that has a mesh size $h$ that is sufficiently small (in our computations eight times smaller).}
We define the $\textcolor{black}{q}$-loss as follows:

\begin{equation}
  \begin{aligned}
    \mathcal{L}(\boldsymbol{\theta}; \mathscr{P}) = \sum_{n=1}^{N} & \bigg[ \norm{\boldsymbol U^n(\boldsymbol{\theta}) - \boldsymbol U_{\textcolor{black}{ref}}^n}_{\textcolor{black}{q}}^{\textcolor{black}{q}} + w_0 \norm{\nabla \boldsymbol U^n(\boldsymbol{\theta}) - \nabla \boldsymbol U_{\textcolor{black}{ref}}^n}_{\textcolor{black}{q}}^{\textcolor{black}{q}} \\
    & + \sum_{\ell_{r,\textcolor{black}{q}} \in \mathcal{R}} w_r \ell_{r,\textcolor{black}{q}}^{\textcolor{black}{q}}(\boldsymbol U^n(\boldsymbol{\theta}), \boldsymbol U_{\textcolor{black}{ref}}^n)\bigg],
  \end{aligned}
\label{eq:loss}
\end{equation}
where $\mathcal{R}$ is a set of regularization terms $\ell_{r,\textcolor{black}{q}}$ in $\textcolor{black}{q}$-norm \textcolor{black}{that will be defined below} and 
\begin{equation}
    \mathscr{P}=\{\mathcal K_h, k, N, \boldsymbol{f}, \Delta t, \boldsymbol U^0, \boldsymbol u_{\textcolor{black}{ref}}\},
    \label{eq:problem-parameters}
\end{equation}
is a tuple that uniquely defines a problem. Namely, a triangulation $\mathcal K_h$, polynomial degree of the basis $k$, number of timesteps $N$, physical flux $\boldsymbol f$, time discretization $\Delta t$, initial condition $\boldsymbol U^0$ and the \textcolor{black}{reference} solution $\boldsymbol u_{\textcolor{black}{ref}}$. Notice that, for the sake of simplicity, we have removed the explicit dependence on $\boldsymbol U^0$ of $\boldsymbol U^n$ for $n > 0$. Indeed, in this loss function, it is implicit the fact that the DG solver described in Section~\ref{sec:algebraic-form} is employed to compute $\boldsymbol U^n$ from $\boldsymbol U^0$. Notice that with a slight abuse of notation, we apply (linear) differential and integral operators to $\boldsymbol{U}^n$, with the meaning of considering the vector of coefficients describing the functions obtained by applying said operator to $\boldsymbol{u}_h(\boldsymbol{x}, t)$. The metrics that we consider \textcolor{black}{(which is a superset of the regularization terms $\ell_{r,\textcolor{black}{q}}$ in $\mathcal{R}$)} are the following:
\begin{itemize}
    \item Error:
    \begin{equation} 
        \epsilon(t^n) = \norm{\boldsymbol U^n - \boldsymbol U_{\textcolor{black}{ref}}^n}_{\textcolor{black}{q}}.
        \label{eq:err-metric}
    \end{equation}
    \item Gradient of the error:
    \begin{equation} 
        \nabla{\epsilon}(t^n) = \norm{\nabla\boldsymbol U^n - \nabla\boldsymbol U_{\textcolor{black}{ref}}^n}_{\textcolor{black}{q}}.
        \label{eq:grad-metric}
    \end{equation}
    \item Jump of the error: 
    \begin{equation} 
        \jmp{\epsilon}(t^n) = \norm{\jmp{\boldsymbol U^n - \boldsymbol U_{\textcolor{black}{ref}}^n}}_{\textcolor{black}{q}}.
        \label{eq:jmp-metric}
    \end{equation}
    The insight of this term comes from the jump-jump term of the DG formulation: we need to penalize large jumps in order to obtain a solution that is not too oscillatory.
    \item Overshoot and undershoot (o/u): 
    \begin{equation} 
      \begin{aligned}
        \text{o/u}(t^n) &= \norm{\boldsymbol U^n \mathbbm{1}_{\boldsymbol U^n > \max \boldsymbol U^n_{\textcolor{black}{ref}}} - \max \boldsymbol U^n_{\textcolor{black}{ref}}}_{\textcolor{black}{q}} \\
        &+ \norm{\boldsymbol U^n \mathbbm{1}_{\boldsymbol U^n < \min \boldsymbol U^n_{\textcolor{black}{ref}}} - \min \boldsymbol U^n_{\textcolor{black}{ref}}}_{\textcolor{black}{q}},
      \end{aligned}   
        \label{eq:ou-metric}
    \end{equation}
    where $\mathbbm{1}_S$ is the indicator function of a set $S$.
    The amplitude of overshoots and undershoots with respect to a reference solution is one of the main criteria with which to evaluate a numerical solution since it provides insights into the behavior and stability of a system. \textcolor{black}{To apply backpropagation to this term, we do not implement it as a product of an indicator function, but we restrict the computations only to the entries of $\boldsymbol{U}^n$ that meet the condition stated in the indicator function.}
    \item Mass variation (mv):
    \begin{equation} 
        \text{mv}(t^n) = \abs{\int_\Omega \boldsymbol U^n - \int_\Omega \boldsymbol U^{n-1}}.
        \label{eq:mv-metric}
    \end{equation}
    We employ this term for problems that have zero mass flux through $\partial \Omega$. Checking that this term is small is key to verifying that the physics of the problem is preserved.
    \item Viscosity penalization (vp):
    \begin{equation} 
        \text{vp}(t^n) = \frac{\frac{1}{\abs{K}}\int_K \mu}{\norm{\nabla \boldsymbol U^n}_{\textcolor{black}{q}} + \varepsilon},
        \label{eq:vp-metric}
    \end{equation}
    where $\varepsilon$ is a small number used for numerical stability, we choose $\varepsilon=10^{-8}$. This term makes sure that we add artificial viscosity only when there are large gradients, therefore encouraging the network to put the least amount possible of viscosity. This term can be substituted by a ``supervised'' term of the kind $|\mu_K - \hat{\mu}_K|$, where $\hat{\mu}_K$ is a reference artificial viscosity given by a classical model. In exchange for a stronger bias on what the NN should learn, this alternative makes the training easier since it steers the learning dynamic towards a known good solution. Still, since this is only a weak constraint, the network is able to discover something new.
\end{itemize}
Not all of them are used in the final version of the loss, indeed we use the minimal number of terms that produce a good model since each term has a weight $w_r$ that needs to be tuned (just once offline). The other terms are still used as metrics to assure the quality of the neural viscosity model. Namely, we have found that the best approach is to employ \textcolor{black}{as the set of regularization terms $\ell_{r,\textcolor{black}{q}}$ in $\mathcal{R}$} the overshoot/undershoot \textcolor{black}{(o/u)} and the viscosity penalization \textcolor{black}{(vp)}. \textcolor{black}{The weights in the loss function for these two terms are hyperparameters to tune.} 

Another key choice is the value of $\textcolor{black}{q}$. We pick $\textcolor{black}{q}=1$ since it tends to emphasize sparsity in the data, meaning it is less sensitive to small values and can be more robust in the presence of outliers. Namely, it is more appropriate in the presence of oscillations that have sparse and isolated peaks, which is our case. \textcolor{black}{Indeed, we have tested other values of $q$, such as $q=2$, and they produced more oscillating results.}

The computations of $\Delta t$ through Eq.~\bref{eq:cfl} must be treated carefully. Indeed, the NN tends to increase to infinity the viscosity $\mu$ in order to have a timestep size $\Delta t$ which is zero. Indeed, as $T=t^N \rightarrow 0$ then the loss tends to zero. There are different solutions to this problem: we can leave the computation outside the AD, we add a penalization term based on $\max_\Omega \mu$ or we can use a fixed $\Delta t$.

\subsection{Training algorithm}\label{sec:training-algo}
\begin{algorithm}[t]
    \caption{Training algorithm}\label{a:train}
    \begin{algorithmic}[1]
    \Procedure{Training}{$\{\mathscr{P}_j\}_{j=1}^{N_p}$, $\{\mathscr{P}_j^\text{test}\}_{j=1}^{N_{pt}}$, \texttt{lr}, scheduler, $N_E$, $N_B$, $\boldsymbol{\theta}$, $\mathcal L$, $L_{\max{}}$, $n_\textnormal{test}$}
        \State $\boldsymbol{\theta^*} \leftarrow \boldsymbol{\theta}, L_{\textnormal{best}} \leftarrow \infty$
        \For{$e$ \textbf{in} $N_E$} \Comment{{\footnotesize Loop over each epoch}} 
            \State \texttt{lr} $\gets$ scheduler(\texttt{lr}) \Comment{Update learning rate according to schedule}
            \State Randomly permute indexes $j$ in $1, ..., N_p$
            \For{$i$ \textbf{in} $N_p / N_B + \textnormal{mod}(N_p, N_B) > 0$} \Comment{Partitions problems in sets of size $N_B$}
                \State Reset gradient accumulation in the optimizer
                \State $L = 0$ \Comment{Reset the loss for the minibatch}
                \For{$j$ \textbf{in} $jN_B, ..., \min\{N_p, (i + 1)N_B\}$}  \Comment{For each problem in the batch}
                        \State $L \leftarrow L + \mathcal L(\boldsymbol{\theta}; \mathscr{P}_{j})$ \Comment{$\mathcal L$ is defined by $q$ and the weights $\{w_i\}_{i \in \mathcal R}$, see Eq.~\bref{eq:loss}}
                \EndFor
                \State $\boldsymbol{\theta} \gets$ AdamW($L$, $\boldsymbol{\theta}$, \texttt{lr}) \Comment{Optimization step }
            \EndFor

            \If{mod$(e, n_\textnormal{test}) = 0$} \Comment{Test outside the automatic differentiation}
                \State $L_\textnormal{test} = \sum_{j=1}^{N_{pt}} \mathcal L(\boldsymbol{\theta}; \mathscr{P}_{j}^\text{test})$
                \If{$L_\textnormal{test} < L_{\textnormal{best}}$}
                    \State $\boldsymbol{\theta^*} \leftarrow \boldsymbol{\theta}$ \Comment{Save model parameters because currently the best}
                \ElsIf{$L_\textnormal{test} > L_{\max{}}$}
                    \State $\boldsymbol{\theta} \leftarrow \boldsymbol{\theta^*}$ \Comment{The model is unstable: load last stable model's parameters}
                    \State \texttt{lr} $\leftarrow$ \texttt{lr}$/ 2$ \Comment{Force learning rate reduction}
                \EndIf
            \EndIf
        \EndFor\\
        \Return $\boldsymbol{\theta^*}$ \Comment{Parameters of the trained NN}
    \EndProcedure
    \end{algorithmic}
\end{algorithm}

\begin{figure}[t]
    \centering
    \includegraphics[width=0.55\textwidth]{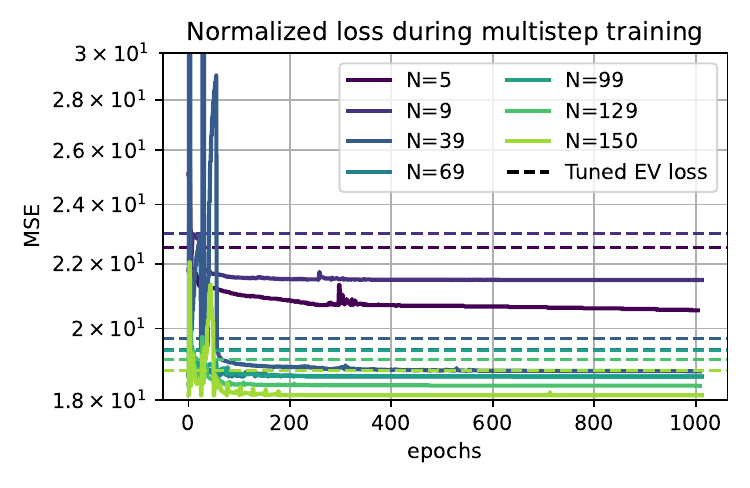}
    \caption{Example of the behavior of the loss for the proposed algorithm during the training for the multi-step procedure described in Section~\ref{sec:training-algo} for a one-dimensional problem.}
    \label{fig:multistep-loss}
\end{figure}

Let us now detail the novel algorithm through which the training is performed. In particular, this training process does not belong to the paradigm of supervised learning since the explicit value of viscosity that the neural network should learn for a certain value of inputs is never computed. \textcolor{black}{Our approach is built on the idea that target viscosity values are not directly used to train the model. Instead, we utilize reference solutions to guide the learning process. However, the approach can neither be classified as unsupervised: when an exact solution is not available, a reference solution must be built using a viscosity model with suitably chosen parameters.} The training algorithm is stated in Algorithm~\ref{a:train} with detailed comments and is summarized in Figure~\ref{fig:train}. It works as follows. It takes as inputs:

\begin{itemize}
    \item A list of $N_p$ tuples of parameters describing the $N_p$ problems $\{\mathscr{P}_j\}_{j=1}^{N_p}$ on which to train the viscosity model on, as defined in Eq.~\bref{eq:problem-parameters}.
    \item A list of $N_{pt}$ problems $\{\mathscr{P}_j^\text{test}\}_{j=1}^{N_{pt}}$ on which to test the algorithm.
    \item An optimizer with its hyperparameters. For instance, we employ AdamW \cite{loshchilov2017decoupled} with learning rate $\texttt{lr}=1\text{e}-3$ and default hyperparameters $\beta_1=0.9, \beta_2=0.999, \varepsilon=10^{-8}, \lambda=0.01$.
    \item A learning rate scheduler \textcolor{black}{\cite{loshchilov2016sgdr}}. We employ a simple ``reduce on plateaus'' \cite{paszke2019pytorch} with patience 30 \textcolor{black}{(number of epochs to wait for a loss reduction)} and reduction factor equal to $\frac{1}{2}$ \textcolor{black}{(factor by which the learning rate is reduced)}.
    \item $N_E$, the number of epochs for which the training lasts.
    \item $N_B$, the batch size of problems considered in one optimization step (gradient computation).
    \item $\boldsymbol{\theta}$, the parameters of the NN. Their shape implicitly defines the architecture of the NN. It is assumed they are already initialized with standard random initialization or a pre-training procedure.
    \item $\mathcal L$, the loss function defined in Eq.~\bref{eq:loss}. It is defined by $\textcolor{black}{q}$ and the weights $\{w_i\}_{i \in \mathcal R}$.
    \item $L_{\max{}}$, the maximum value of the loss. If the loss is larger than this threshold the model is considered unstable. We use $L_{\max{}}=10^{30}$.
    \item $n_{\textnormal{test}}$, the frequency with which the model is tested.
\end{itemize}
Then, the training algorithm closely resembles a supervised learning loop, where the loss is substituted with the one we described in the previous section. Namely, we create minibatches of parameters and make updates based on these parameters. The key point being that we do not have a target value that the NN should learn but a loss that describes the physics of the problem through a DG solver. Indeed, the peculiarity of this training algorithm is that, like in RL, it learns from the consequences of the actions (choices of $\mu_K$) it makes in the DG solver, creating a feedback loop. In particular, the value of artificial viscosity at the timestep $n$ influences the solution at $n+1$ which in turn is used as input of the network to determine the next value of the viscosity. This enables the model to learn long-term dependencies and optimal strategies in response to evolving scenarios. 

Among the most interesting features of this algorithm is that, compared to the supervised training process, we can avoid designing and building a dataset. This is particularly expensive from the point of view of human time since it requires tuning the artificial viscosity parameters by hand and picking the best choice by confronting the various solutions. In particular, differently from the classical loss Eq.~\bref{eq:loss_classic}, our loss Eq.~\bref{eq:loss} necessitates only a reference solution, which is much cheaper to build. Coincidentally, this feature also solves the aforementioned drawback (C), since the neural network automatically learns the optimal viscosity values. Since the problem is highly non-convex, this cannot be a global optimum, however, we guarantee that this is at least a local minimum.

Let us now briefly discuss the choice of hyperparameters for this training algorithm. We observe that we obtain better results with a small batch size, usually between two and four. It is paramount to choose $N$ properly: if it is too small, the neural network does not have enough data to generalize, if it is too large we may incur into exploding gradient problems unless the neural network has a proper weight initialization. We propose a multi-step pre-train in which we iterate the training Algorithm~\ref{a:train} for increasing values of $N$ for a few epochs, starting from $N=N_0$ up to the desired final $N$ with a certain step $N_e$. The increase in number of timesteps $N_e$ is computed adaptively, namely\textcolor{black}{, defined}
\begin{equation*}
    \textcolor{black}{\texttt{clip}(a, b, c) = \max(\min(b, c), a),}
\end{equation*}
\textcolor{black}{we compute $N_e$ as} $N_e = \texttt{clip}(N_e^{\min{}}, N_e^\alpha, N_e^{\max{}})$, where $N_e^\alpha$ is the largest integer so that the average loss per timestep with $N + N_e^\alpha$ timesteps is less than a constant (hyperparameter) $\alpha_e > 0$ multiplied by the average loss per timestep with $N$ timesteps. Once the NN reaches good generalization properties it is observed that the average loss per timestep decreases when $N$ increases. Indeed, this shows that the neural viscosity model is making accurate predictions of the ``future'' (data it has never seen before). Figure~\ref{fig:multistep-loss} shows an example of the training loss. This process is expensive from a computational point of view, however, it enables the NN to discover a viscosity model without any biases.

Since the magnitude of the loss changes significantly between the start and the end of the optimization process, a key ingredient to increase the stability and efficiency of the algorithm is a learning rate scheduler. A simple scheduler that geometrically reduces the error on plateaus worked best for us, even compared to more modern and popular alternatives like the cosine scheduler \cite{loshchilov2016sgdr}. This also allows us to have an algorithm that is less influenced by the starting value of the learning rate since, if it is too large, is automatically decreased. In our experiments, we employed a starting learning rate of $5\cdot10^{-2}$.

Finally, we stress that, despite all these precautions, it seldom happens that a gradient update produces a viscosity model that is not able to handle the Gibbs phenomenon. For this reason, it is important to constantly monitor the loss, save the best parameters of the model, and use them for a restart with a reduced learning rate. This simple but effective procedure completely eradicates this issue.

\subsection{Training in practice}\label{sec:training-in-practice}
\begin{table}[t]
    \centering
    \begin{minipage}{0.5\textwidth}
        \centering
        \footnotesize
        \setlength\extrarowheight{3pt}
        \begin{tabular}{l}
            \hline 
            \textbf{1D initial conditions}\\
            \hline 
            $u_0(x) = \frac{\omega}{2} \sin(\omega \pi x), \quad \omega = 1, 3, 5$\\
            $u_0(x) = \mathbbm{1}_{[\frac{1}{4}, \frac{3}{4}]}$\\
            $u_0(x) = 10 (\frac{1}{2} - |x - \frac{1}{2}|)$\\
            $u_0(x) = e^{-100(x - \frac{1}{2})^2}$\\
            $u_0(x) = \omega_1 \mathbbm{1}_{[0, \frac{1}{5}]} + \omega_2 \mathbbm{1}_{[\frac{1}{5}, \frac{2}{5}]} + \omega_3 \mathbbm{1}_{[\frac{3}{5}, 1]}, \quad \omega_i = -4, 6, 10$\\
            $u_0(x) = \sin(\omega\pi x) \mathbbm{1}_{[\frac{1}{4}, \frac{1}{2}]} + \sin(2\omega\pi x) \mathbbm{1}_{[\frac{1}{2}, \frac{3}{4}]}, \quad \omega = 4, 8$\\
            $u_0(x) = -\sin(6 \pi x) \mathbbm{1}_{[\frac{1}{6}, \frac{5}{6}]}$\\
            \makecell[l]{$u_0(x) = \omega_1(x - \frac{1}{6}) \mathbbm{1}_{[\frac{1}{6}, \frac{1}{3}]} +  \omega_2(x - \frac{1}{2})\mathbbm{1}_{[\frac{1}{3}, \frac{2}{3}]} $\\
            $\qquad\quad +  \omega_3(x - \frac{5}{6})\mathbbm{1}_{[\frac{2}{3}, \frac{5}{6}]}, \quad \omega_i = 2, 6, 10$}\\
            $u_0(x) = (16\abs{x - \frac{1}{2}} - 2) \mathbbm{1}_{[\frac{1}{4}, \frac{3}{4}]}$\\
            \hline 
        \end{tabular}
    \end{minipage}%
    \begin{minipage}{0.5\textwidth}
        \centering
        \footnotesize
        \setlength\extrarowheight{3pt}
        \begin{tabular}{l}
            \hline 
            \textbf{2D initial conditions}\\
            \hline 
            $u_0(\boldsymbol x) = \frac{\omega_1\omega_2}{2} \sin(\omega_1 \pi x_1)\sin(\omega_2 \pi x_2), \quad \omega_i = 1, 3, 5$\\
            $u_0(\boldsymbol x) = \mathbbm{1}_{[\frac{1}{4}, \frac{3}{4}]}(x_1)\mathbbm{1}_{[\frac{1}{4}, \frac{3}{4}]}(x_2)$\\
            $u_0(\boldsymbol x) = e^{-100\norm{\boldsymbol x - \frac{1}{2} \boldsymbol{1}}^2}$\\
            \makecell[l]{$u_0(\boldsymbol x) = \omega_1 \mathbbm{1}_{[0, \frac{1}{5}]}(x_1) \mathbbm{1}_{[\frac{3}{5}, 1]}(x_2) $\\
            $\qquad\quad+ \omega_2 \mathbbm{1}_{[\frac{1}{5}, \frac{2}{5}]}(x_1)\mathbbm{1}_{[\frac{1}{5}, \frac{2}{5}]}(x_2) $\\
            $\qquad\quad+ \omega_3 \mathbbm{1}_{[\frac{3}{5}, 1]}(x_1) \mathbbm{1}_{[0, \frac{1}{5}]}(x_2), \quad \omega_i = -4, 6, 10$}\\
            \makecell[l]{$u_0(\boldsymbol x) = \sin(\omega\pi x) \mathbbm{1}_{[\frac{1}{4}, \frac{1}{2}]}(x_1)\mathbbm{1}_{[\frac{1}{4}, \frac{1}{2}]}(x_2) $\\
            $\qquad\quad + \sin(2\omega\pi x) \mathbbm{1}_{[\frac{1}{2}, \frac{3}{4}]}(x_1)\mathbbm{1}_{[\frac{1}{2}, \frac{3}{4}]}(x_2), \quad \omega = \frac{5}{2}, 4, 8$}\\
            \makecell[l]{$u_0(\boldsymbol x) = (\omega_1 x_1 + \omega_2 x_2 - \frac{1}{4}) \mathbbm{1}_{[\frac{1}{4}, \frac{3}{4}]}(x_1) \mathbbm{1}_{[\frac{1}{4}, \frac{3}{4}]}(x_2)$\\
            $\qquad\quad \omega_i = -1, 0, 1, 4$}\\
            \hline 
        \end{tabular}
    \end{minipage}%
    \caption{Initial conditions $u_0$ employed in the training phase for advection and Burgers' fluxes.}\label{tab:u0}
\end{table}

\begin{figure}[t]
    \centering
    \includegraphics[width=0.6\textwidth]{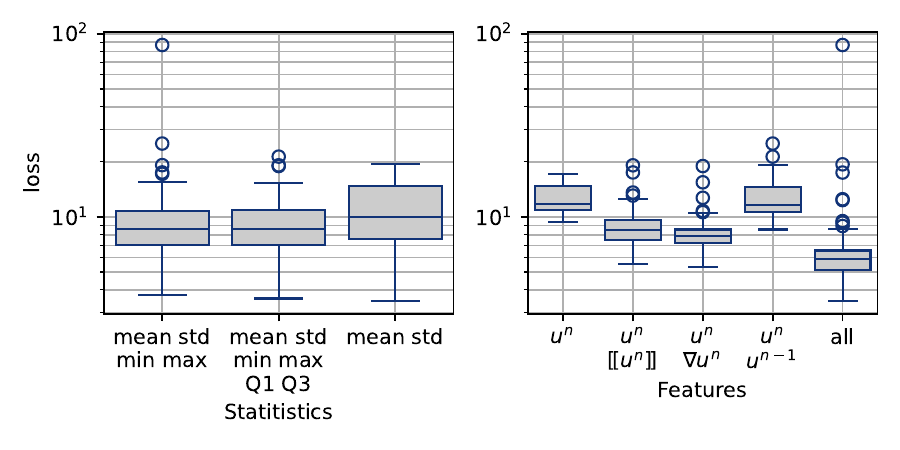}
    \caption{Loss of the trained viscosity model depending on the nodal quantities used as input of the NN (\textit{right}) and on the functions used to agglomerate those values (\textit{left}). See Section~\ref{sec:neural-viscosity} for a complete definition of these quantities. The results are obtained after one step of the random search relative to the NN employed in the 1D cases.}
    \label{fig:hyperparameters-rand-search-pt0}
\end{figure}

\begin{figure}[t]
    \centering
    \includegraphics[width=0.75\textwidth]{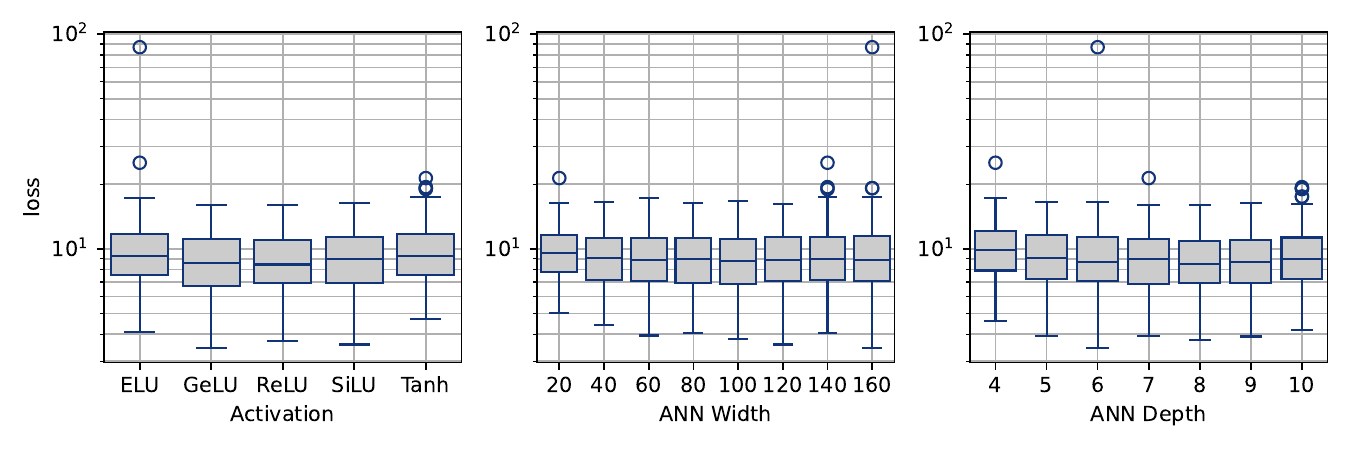}
    \caption{Loss of the trained viscosity model depending on hyperparameters defining the architecture of the NN agglomerated in boxplots. The results are obtained after one step of the random search relative to the NN employed in the one-dimensional cases.}
    \label{fig:hyperparameters-rand-search-pt1}
\end{figure}

In this section, we provide more details on how the training is carried out in the one-dimensional ($d=1$) case. This procedure can be easily generalized to different dimensions. The first step is to choose which problems to train the neural network on, namely, we need to choose different combinations 
$\{\mathscr{P}_j\}_{j=1}^{N_p}, \{\mathscr{P}_j^\text{test}\}_{j=1}^{N_{pt}}$ of problems on which to train and test. 
The two sets are an 80-20 random split of a union of cartesian products of parameter spaces. For the 1D case we consider uniform discretizations of $(0, 1)$ with $h=\frac{2^{-i}}{10}, i=1,...,4$, degrees $k=1,2,4$, the flux $\boldsymbol f$ of linear advection, Burgers' equation and Euler system, a constant timestep $\Delta t$ that is the largest such that the considered problem is stable (it varies from problem to problem), number of timesteps $N=100$ for linear advection and appositely chosen as the smallest $N\in\mathbb{N}$ such that shocks and waves appear in Burgers' and Euler cases. Concerning initial conditions, we use the functions reported on the left of Table~\ref{tab:u0} for advection and Burgers' problems, and a Sod shock tube problem \cite{sod1978survey} for Euler's. A representation of the initial conditions is shown in Figure~\ref{fig:ic-adv-1d-train}. For the 2D case, the problems are chosen in an analogous way apart from the following parameters. Since unstructured grids offer flexibility in discretizing computational domains, allowing for efficient representation of intricate shapes and unbounded regions, we use both structured and unstructured triangulations $\mathcal K_h$ to train the NN. Namely, we discretize $(0, 1)^2$ with $h=\frac{2^{-i}}{10}, i=1,...,3$. The initial conditions are shown on the right of Table~\ref{tab:u0} for advection and Burgers' equations. We also add the Euler's system in the Riemann configuration 12 of \cite{kurganov2002solution} to the training environment.

Once the parameters defining the training problems are chosen, we pass to hyperparameter tuning. As it is commonly done, we employ a random search: we specify the range or distribution for each hyperparameter we want to tune, and we define the search space as the cartesian product of these ranges. The model is trained using random choices of hyperparameter values, and the loss is used as an evaluation metric. This process is repeated while, at each iteration, reducing the search space and increasing the fraction of the combinatorial space probed and the number of epochs $N_E$.

As an optimizer, we employ AdamW (Adam with Weight Decay) \cite{loshchilov2017decoupled}. In our experiments, it significantly outperformed stochastic gradient descent (SDG) and proved to be less sensitive to the choice of the learning rate with respect to Adam.

We then focus on five hyperparameters that we consider to be particularly relevant, namely the width and depth of the NN, the activation function, the local features extracted from each element $K$, and which statistics are extracted from these features. In Figure~\ref{fig:hyperparameters-rand-search-pt0}~and~\ref{fig:hyperparameters-rand-search-pt1} we show the results after one iteration of random search (about ten thousand combinations were tested). After three iterations of random search we, conclude the following (for the one-dimensional case): 
\begin{itemize}
    \item using the mean, standard deviation, minimum, and maximum of the nodal values as describing statistics provides the best balance between computational cost and accuracy;
    \item the gradient of the \textcolor{black}{state variable}, the jump across the element, and the \textcolor{black}{state variable} at the previous timestep are all features that significantly help the training;
    \item the Gaussian error linear unit (GeLU) \cite{hendrycks2016gaussian} is the activation function that \textcolor{black}{we choose since it reached the smallest loss} when compared to the rectified linear unit (ReLU), exponential linear unit (ELU), sigmoid linear unit (SiLU) and hyperbolic tangent (Tanh);
    \item a large model, with a width greater than 100, \textcolor{black}{has a smaller loss on average}. In particular, we choose to use 120 neurons per layer;
    \item models with depth between six and eight (which are neither shallow nor very deep) reach the lowest loss on average. In our experiments, we employ seven layers.
\end{itemize}

\textcolor{black}{Therefore, the input vector of the neural network has size 19 (4 for the reference variable, its gradient, its value at the previous time step, its flux; 2 for the jump; 1 for the polynomial degree) in 1D and 29 in 2D (4 for the reference variable, its gradient, its value at the previous time step, its flux and its jump; 1 for the polynomial degree)}

\section{Numerical Results}\label{sec:results}
\begin{figure}[t]
    \centering
    \begin{minipage}{0.5\textwidth}
        \centering
        \includegraphics[width=0.8\linewidth,trim={0 0 0 .7cm},clip]{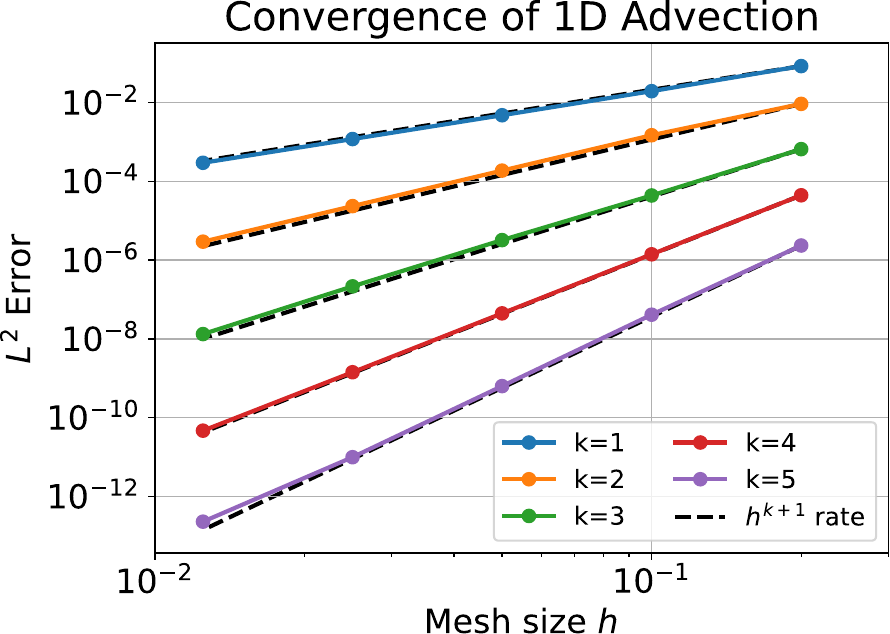}
        \caption{Computed errors and convergence rates for test case of linear advection ($d=1$).}
        \label{fig:convergence-1d-adv}
    \end{minipage}%
    \begin{minipage}{0.5\textwidth}
        \centering
        \includegraphics[width=0.8\linewidth,trim={0 0 0 .7cm},clip]{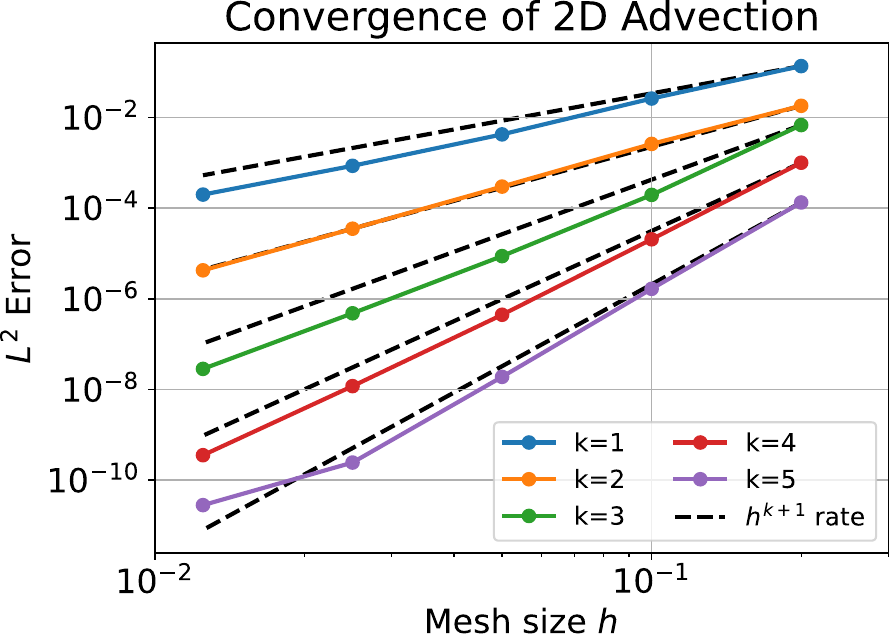}
        \caption{Error estimates and convergence rates for the 2D smooth test case of linear advection.}
        \label{fig:convergence-2d-adv}
    \end{minipage}
\end{figure}

In this section, our objective is to assess the accuracy of the proposed training algorithm in practice. In the following, we present a comprehensive set of numerical results, encompassing both scalar equations and systems, to showcase the efficacy of physics-informed machine learning. All the tests have been implemented in an appositely developed Python library for the solution
of conservation laws with the DG method that employs Pytorch \cite{paszke2019pytorch} for automatic differentiation and the implementation of the NNs.

First, we present a convergence test to verify the functioning of our numerical solver and neural viscosity model. Then, we present tests carried out on problems that are not included in the training environment and are characterized by different parameters, such as: flux, initial condition, boundary conditions, domain discretization, polynomial degree, and CFL. In this way, we can assess the generalization capabilities of the trained networks and test the robustness of our algorithm in an impartial way. Our findings consistently demonstrate performance superior to the optimally tuned EV models. This is especially relevant considering that tuning is itself a significantly expensive operation.

\subsection{Verification tests: convergence rate}

Here, we aim to verify the accuracy of the presented approach. In particular, we assess the ability of the neural viscosity model to maintain the expected accuracy for a smooth solution. Indeed, the addition of a dissipative term might adversely affect the accuracy. The assessment of the numerical scheme employs the $L^2(\Omega)$-norm of the discretization error vector, evaluated at the final time $T$. Namely, theoretical results for inviscid hyperbolic problems \cite{hesthaven2007nodal,hesthaven2017numerical} show that the error
\begin{equation*}
    \epsilon = \norm{u_h (\cdot, T ) - u(\cdot, T )}_{L^2(\Omega)} \lesssim h^{k+1},
\end{equation*}
where the hidden constant depends on $k$, given that $u$ is sufficiently regular and a Rusanov flux is employed.

\subsubsection{Test case 1: one-dimensional smooth problem}
To test the convergence we consider Problem~\bref{eq:strong_conservation_law} with the linear flux defined by Eq.~\bref{eq:adv-1d}, namely with a constant unitary transport field and with the smooth initial conditions
\begin{equation*}
    u_0(x) = \frac{1}{2} + \sin(2\pi x).
\end{equation*}
The problem is endowed with periodic boundary conditions.
Concerning the discretization, we use the spatial domain $\Omega = (0, 1)$, $\text{CFL}=0.05$, final time $T=0.4$ and Runge-Kutta of the fourth-order. In Figure~\ref{fig:convergence-1d-adv} we report the computed errors in the $L^2(\Omega)$-norm, together with the expected rates of convergence as a function of the mesh size $h$ and the polynomial degree $k$. We can notice that we recover the expected convergence rate.

\subsubsection{Test case 2: two-dimensional smooth problem}
Similarly, we extend the previous test case in two dimensions. Namely, we consider Problem~\bref{eq:strong_conservation_law} with constant unitary transport field $\boldsymbol \beta = [1, 1]^\top$ as defined in Eq.~\bref{eq:advection-2d} with smooth initial conditions
\begin{equation*}
    u_0(x_1, x_2) = \frac{1}{2} + \sin(2\pi x_1)\sin(2\pi x_2).
\end{equation*}
The problem is endowed with periodic boundary conditions.
Concerning the discretization, we employ a structured triangular grid of the domain $\Omega = (0, 1)^2$ of granularity $h=\frac{\sqrt{2}}{2^i10}, i = 1, ..., 5$. In Figure~\ref{fig:convergence-2d-adv} we report the computed errors in the $L^2(\Omega)$-norm and the expected rates of convergence. The NN technique achieves the optimal convergence rate. Numerical experiments show also that we obtain the same \textcolor{black}{rate of convergence} independently of the mesh type (structured or unstructured).

\subsection{One-dimensional problems}
\begin{figure}[t]
    \centering
    \begin{minipage}{0.5\textwidth}
        \centering
        \includegraphics[width=0.9\linewidth]{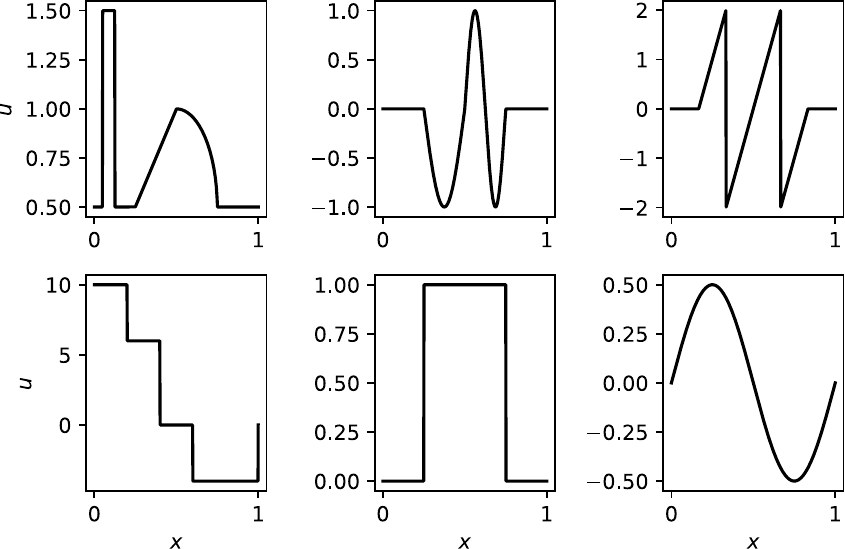}
        \caption{Examples of initial conditions $u_0$ used to train the NN.}
        \label{fig:ic-adv-1d-train}
    \end{minipage}%
    \begin{minipage}{0.5\textwidth}
        \centering
        \includegraphics[width=0.8\linewidth]{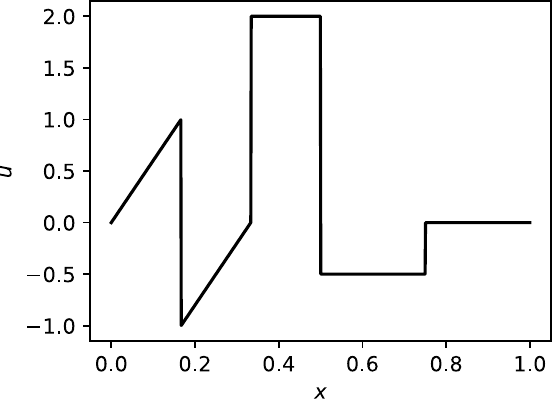}
        \caption{Initial condition $u_0$ for test case 3 and 4.}
        \label{fig:ic-adv-1d}
    \end{minipage}
\end{figure}

To gain more insight into how the algorithm works, we first analyze one-dimensional problems. The following test cases can be divided into two categories. The first two deal with scalar problems. To test the generalization properties of the NN, we consider an initial condition, grid spacing, polynomial degree, and CFL not included in the training environment. In particular, we consider Problem~\bref{eq:strong_conservation_law} with flux defined by Eq.~\bref{eq:adv-1d} and Eq.~\bref{eq:burgers-1d} and we employ as initial condition the function
\begin{equation}
    u_0(x) = 
    \begin{cases}
        6x & x \in (0, \frac{1}{6}],\\
        6(x - \frac{1}{3}) & x \in (\frac{1}{6}, \frac{1}{3}],\\
        2 & x \in (\frac{1}{3}, \frac{1}{2}],\\
        -\frac{1}{2} & x \in (\frac{1}{2}, \frac{3}{4}],\\
        0 & x \in (\frac{3}{4}, 1),\\
    \end{cases}
    \label{eq:ic-t34}
\end{equation}
which is represented in Figure~\ref{fig:ic-adv-1d}. Instead, the other two test cases seek validation in a more challenging setting. Namely, we consider two common tests for the accuracy of the solution of the Euler equations: the Sod shock tube problem \cite{sod1978survey} and the Shu-Osher problem \cite{shu1988efficient}.

\subsubsection{Test case 3: linear advection}
\begin{figure}[!t]
    \centering
    \begin{minipage}{0.33\textwidth}
        \centering
        \includegraphics[width=1.0\linewidth]{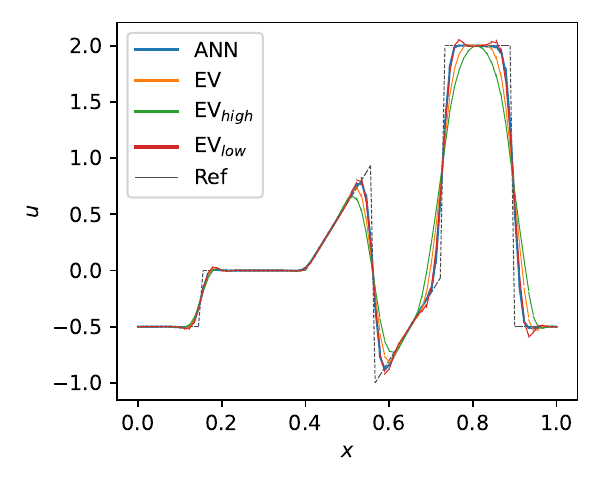}
    \end{minipage}%
    \begin{minipage}{0.33\textwidth}
        \centering
        \includegraphics[width=1.0\linewidth]{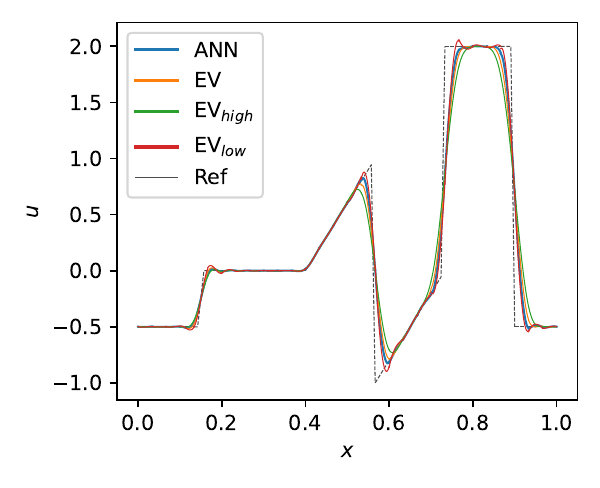}
    \end{minipage}%
    \begin{minipage}{0.33\textwidth}
        \centering
        \includegraphics[width=1.0\linewidth]{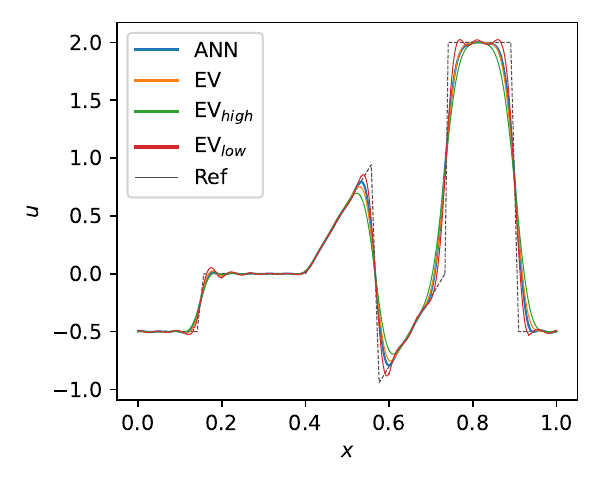}

    \end{minipage}
    \caption{Test case 3: solution at $T=0.4$ for $k=1,3,5$ (from left to right).}
    \label{fig:t3-ufin}
\end{figure}

\begin{figure}[!t]
    \centering
    \begin{minipage}{0.36\textwidth}
        \centering
        \includegraphics[height=7.5cm]{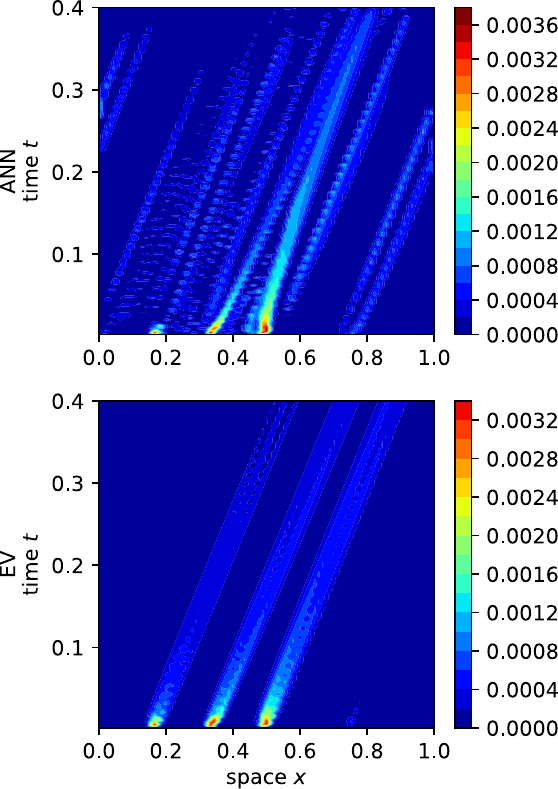}
    \end{minipage}%
    \begin{minipage}{0.32\textwidth}
        \centering
        \includegraphics[height=7.5cm]{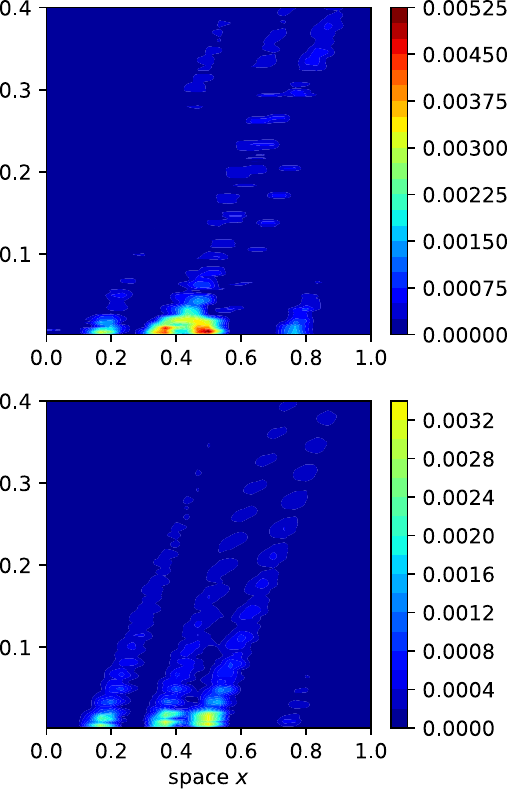}
    \end{minipage}%
    \begin{minipage}{0.32\textwidth}
        \centering
        \includegraphics[height=7.5cm]{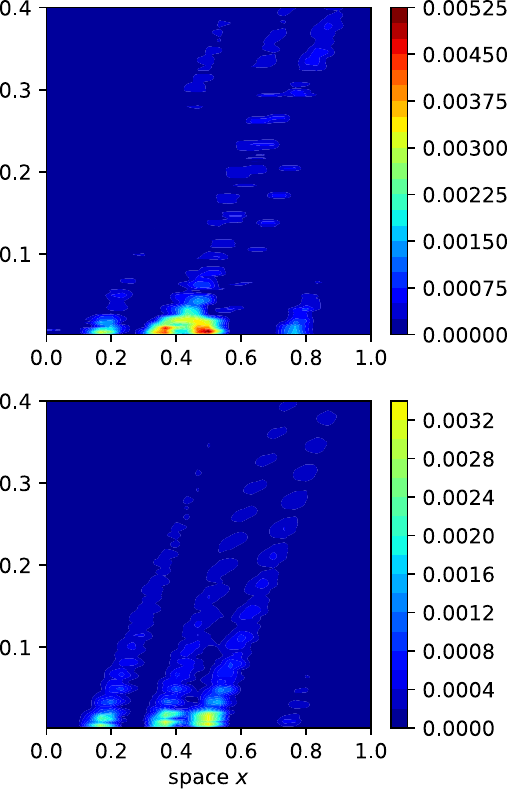}
    \end{minipage}
    \caption{Test case 3: comparison of space-time artificial viscosity for $k=1,3,5$ (from left to right).}
    \label{fig:t3-mu}
\end{figure}


\begin{table}[t]
    \centering
    \footnotesize
    \begin{tabular}{l|cc|cc|cc}
                        & \multicolumn{2}{c}{$k=1$} & \multicolumn{2}{c}{$k=3$} & \multicolumn{2}{c}{$k=5$} \\
                        & NN         & EV        & NN         & EV        & NN         & EV        \\
        \hline
        $\epsilon$       & 2.2311e+03 & 3.1745e+03 & 1.7990e+03 & 2.1943e+03 & 1.5680e+03 & 1.6753e+03 \\
        $\nabla\epsilon$ & 1.7306e+03 & 1.9666e+03 & 6.7370e+03 & 6.9916e+03 & 1.1302e+04 & 1.1380e+04 \\
        $\jmp{\epsilon}$ & 1.5933e+02 & 1.1955e+02 & 4.2461e+01 & 3.4186e+01 & 2.0260e+01 & 2.0562e+01 \\
        o/u              & 3.4889e+00 & 1.2134e+01 & 2.2778e+00 & 6.2186e+00 & 2.5903e+00 & 7.3782e+00 \\
        mv               & 4.0360e-03 & 4.5267e-03 & 1.6535e-03 & 2.8121e-03 & 1.4988e-03 & 1.6169e-03
    \end{tabular}
    \caption{Test case 3: comparison of cumulative $L^1$ error metrics Eq.~\bref{eq:err-metric}~-~\bref{eq:mv-metric} over all timesteps from $0$ to $T$.}
    \label{tab:t3-err}
\end{table}

The first test case we consider is a linear advection problem with a highly discontinuous initial condition, namely \bref{eq:ic-t34} in the $\Omega = (0, 1)$ domain. The simulation is performed until $T=0.4$ with polynomial degrees $k=1,3,5$ and endowed with periodic boundary conditions. To compare solutions with a similar number of degrees of freedom we select $h=\frac{1}{60}, \frac{1}{30}, \frac{1}{15}$, $\text{CFL}=0.2,0.5,0.75$ \textcolor{black}{and $C=0.2,0.0556,0.03$}, respectively. The solution obtained with the neural viscosity model is compared with a tuned EV model ($c_K=0.6$, $c_{\max{}}=0.3$, see Eq.~\bref{eq:max_visc},~\bref{eq:ev-visc}) and two instances of EV model that are not optimally tuned, namely one case where the injected dissipation is too small and too large. In Figure~\ref{fig:t3-ufin} we compare the solutions at the final time: we can see that the NN based solution \textcolor{black}{preserves the same symmetry properties of the exact solution}, essentially non-oscillatory, and well captures the front of the wave. By analyzing the space-time plot of the viscosity introduced via the NN, as shown in Figure~\ref{fig:t3-mu}, we notice that the NN tends to add a larger amount of \textcolor{black}{less} localized viscosity when compared with EV \textcolor{black}{for $k=1$}. \textcolor{black}{It is instead more similarly localized for $k=3,5$.} More quantitative results about the error are shown in Table~\ref{tab:t3-err}, we notice that the NN always outperforms the EV model in terms of $L^1$ error and overshoots/undershoots. It also obtains comparable results when considering mass variation (mv), an important metric to asses that the physics of the problem is respected. \textcolor{black}{This comparison is limited by the choice of parameters $c_K$ and $c_{\max{}}$ we made and by the considered metrics. Indeed, while we manually tuned these parameters, there might exist a choice that leads to smaller error metrics for the EV model.}

\subsubsection{Test case 4: Burgers' equation}

\begin{figure}[!t]
    \centering
    \begin{minipage}{0.33\textwidth}
        \centering
        \includegraphics[width=1.0\linewidth]{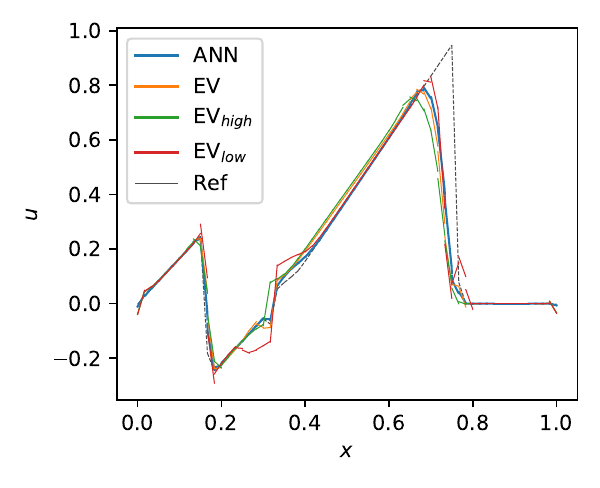}
    \end{minipage}%
    \begin{minipage}{0.33\textwidth}
        \centering
        \includegraphics[width=1.0\linewidth]{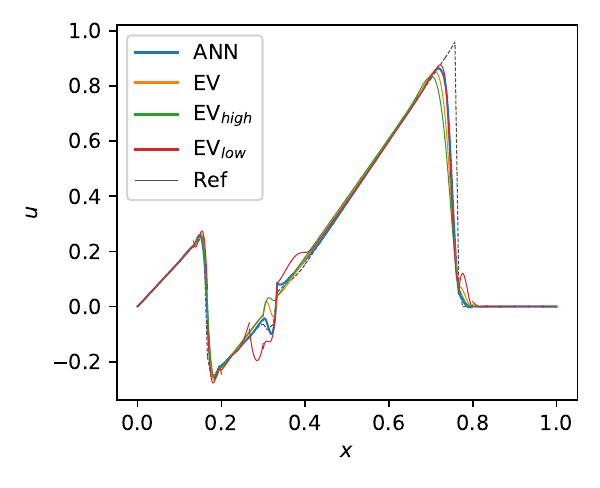}
    \end{minipage}%
    \begin{minipage}{0.33\textwidth}
        \centering
        \includegraphics[width=1.0\linewidth]{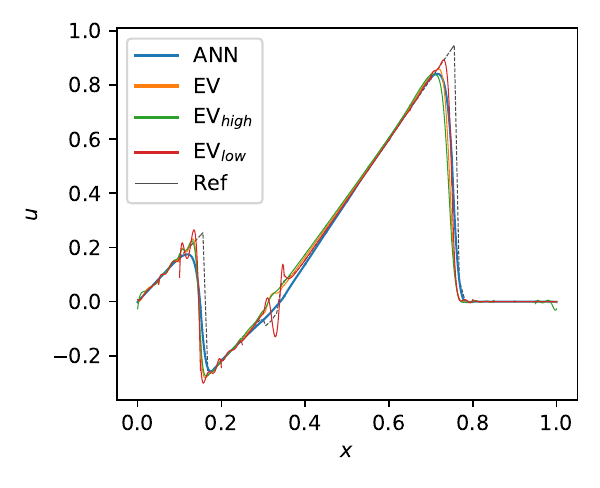}
    \end{minipage}
    \caption{Test case 4: solution at final time for $k=1,3,5$ (from left to right).}
    \label{fig:t4-ufin}
\end{figure}

\begin{figure}[!t]
    \centering
    \begin{minipage}{0.36\textwidth}
        \centering
        \includegraphics[height=7.5cm]{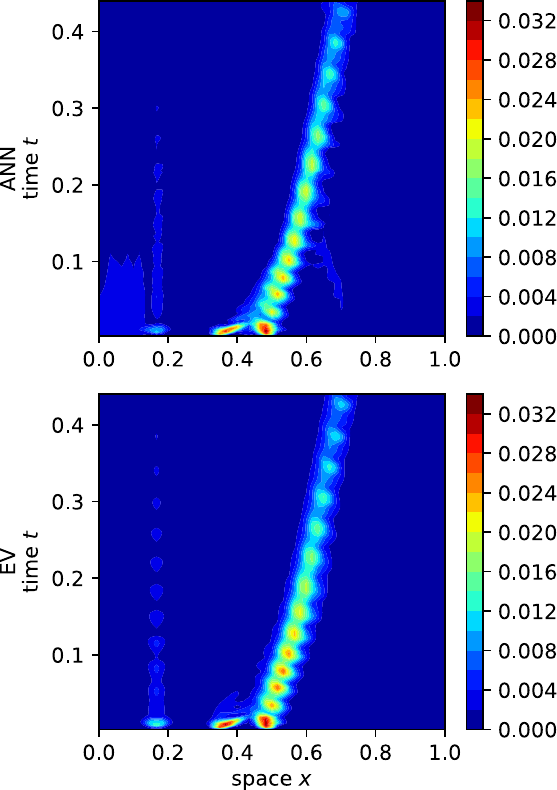}
    \end{minipage}%
    \begin{minipage}{0.32\textwidth}
        \centering
        \includegraphics[height=7.5cm]{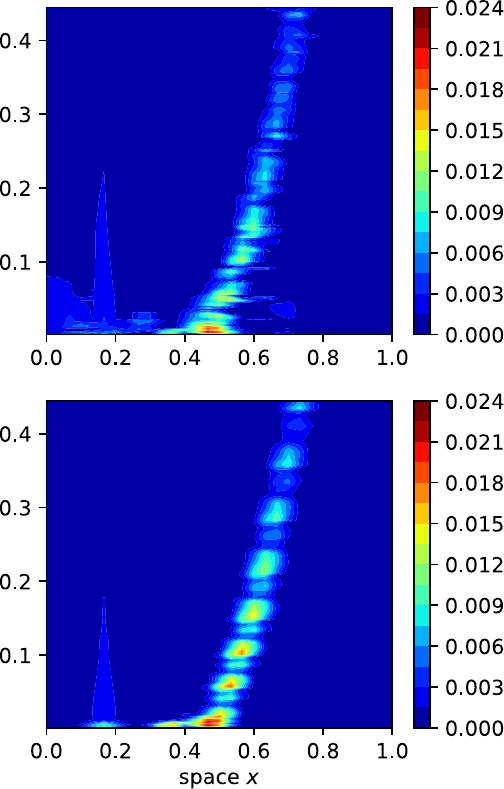}
    \end{minipage}%
    \begin{minipage}{0.32\textwidth}
        \centering
        \includegraphics[height=7.5cm]{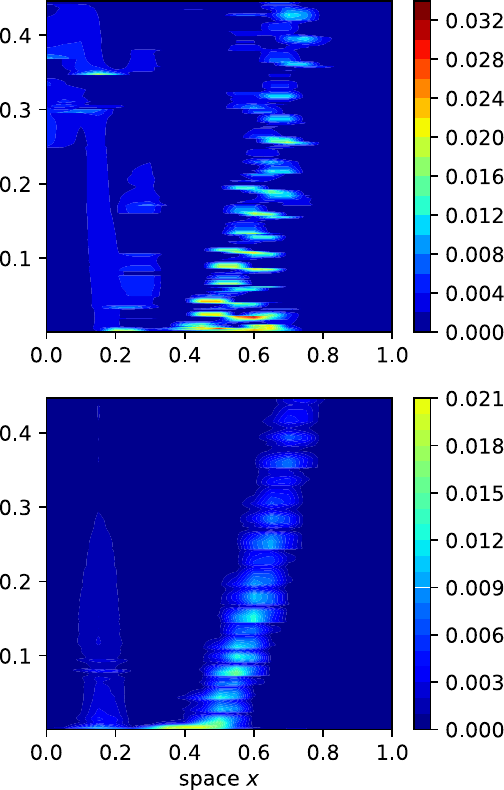}
    \end{minipage}
    \caption{Test case 4: comparison of space-time artificial viscosity for $k=1,3,5$ (from left to right).}
    \label{fig:t4-mu}
\end{figure}


\begin{table}[!t]
    \centering
    \footnotesize
    \begin{tabular}{l|cc|cc|cc}
                        & \multicolumn{2}{c}{$k=1$} & \multicolumn{2}{c}{$k=3$} & \multicolumn{2}{c}{$k=5$} \\
                        & NN         & EV        & NN         & EV        & NN         & EV        \\
        \hline
        $\epsilon$       & 8.7204e+02 & 1.1423e+03 & 9.7797e+02 & 1.3587e+03 & 2.6155e+03 & 3.2874e+03\\
        $\nabla\epsilon$ & 5.3064e+02 & 6.1239e+02 & 3.1401e+03 & 3.5330e+03 & 1.4578e+04 & 1.7510e+04\\
        $\jmp{\epsilon}$ & 1.0947e+02 & 1.1138e+02 & 4.4744e+01 & 5.0420e+01 & 3.4530e+01 & 1.2684e+02\\
        o/u              & 6.5979e-01 & 9.8408e-01 & 6.5894e-01 & 1.3003e+00 & 9.7566e-01 & 2.5800e+00\\
    \end{tabular}
    \caption{Test case 4: comparison of cumulative $L^1$ error metrics Eq.~\bref{eq:err-metric}~-~\bref{eq:ou-metric} over all timesteps from $0$ to $T$.}
    \label{tab:t4-err}
\end{table}

In this test case, we keep the same parameters of the previous test case but we change the physical flux, namely, we employ a quadratic one, and the $\text{CFL}=0.15,0.4,0.4$. \textcolor{black}{The Courant number is $C=0.3317,0.0924,0.0302$.} These settings are more challenging due to the fact that as the solution evolves rarefaction fans and shock waves develop and collide. We compare the NN solution with a tuned EV model ($c_K=3.0, c_{\max{}}=1.0$) and two other not optimally tuned EV models. Figure~\ref{fig:t4-ufin} shows that one of the EV models ($c_K=1.0, c_{\max{}}=0.5$) is not able to add enough viscosity thus leading to large oscillation. The parameters used in this case are not far from the manually tuned: this showcases how time-consuming the parameter selection can be.
In this case, the NN model is by far the best in avoiding overshoots/undershoots while correctly predicting the profile of the wave. The history of the viscosity in Figure~\ref{fig:t4-mu} shows that as in the previous test cases, the NN tends to inject more viscosity than the EV model, \textcolor{black}{in a more spread} and less continuous way. For $k=5$ in particular, it is possible to see that the NN also tends to add dissipation in $x=0.2, t=0.35$ where the EV model did not. It remains to be investigated if this is a generalization error (since $k=5$ was not present in the training problems) or a feature of the discovered neural viscosity model. As shown in Table~\ref{tab:t4-err}, the NN model outperforms all the EV models under all the considered metrics. This is an especially interesting result when considering that the reduction of the $L^1$ error for $k=3$ is of 25\%, while also reducing overshoots/undershoots of almost 50\%.

\subsubsection{Test case 5: Sod shock tube}

\begin{figure}[t]
    \centering
    \begin{minipage}{0.33\textwidth}
        \centering
        \includegraphics[width=0.95\linewidth]{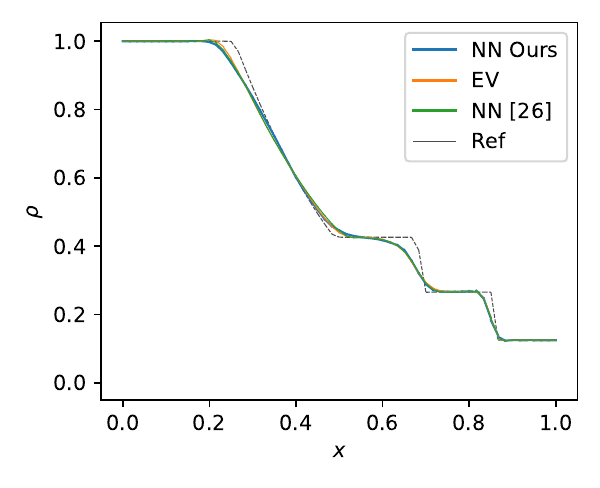}
    \end{minipage}%
    \begin{minipage}{0.33\textwidth}
        \centering
        \includegraphics[width=0.95\linewidth]{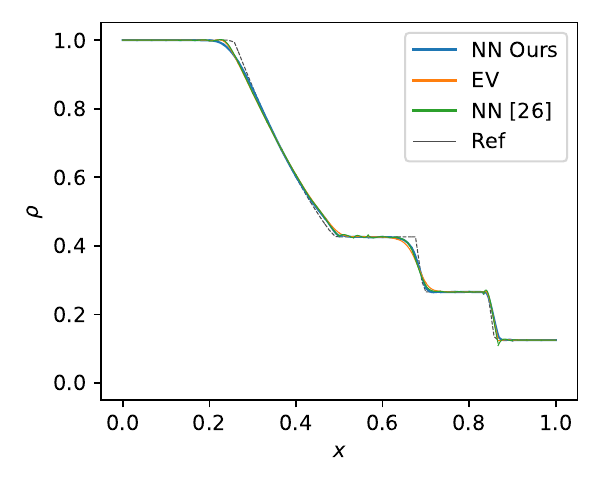}
    \end{minipage}%
    \begin{minipage}{0.33\textwidth}
        \centering
        \includegraphics[width=0.95\linewidth]{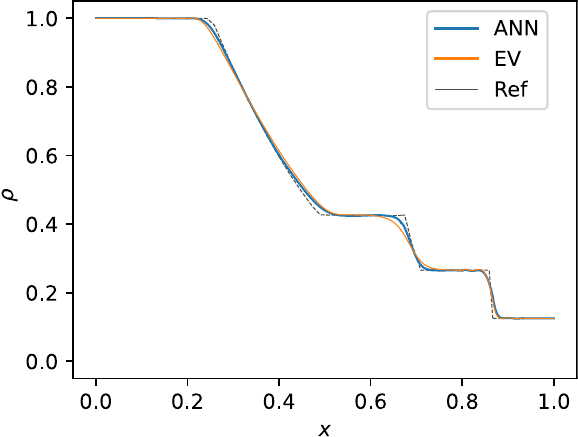}
    \end{minipage}
    \caption{Test case 5: density $\rho$ at $T=0.2$ for $k=1,3,5$ (from left to right).}
    \label{fig:t5-ufin}
\end{figure}

\begin{figure}[t]
    \centering
    \begin{minipage}{0.36\textwidth}
        \centering
        \includegraphics[height=11.25cm]{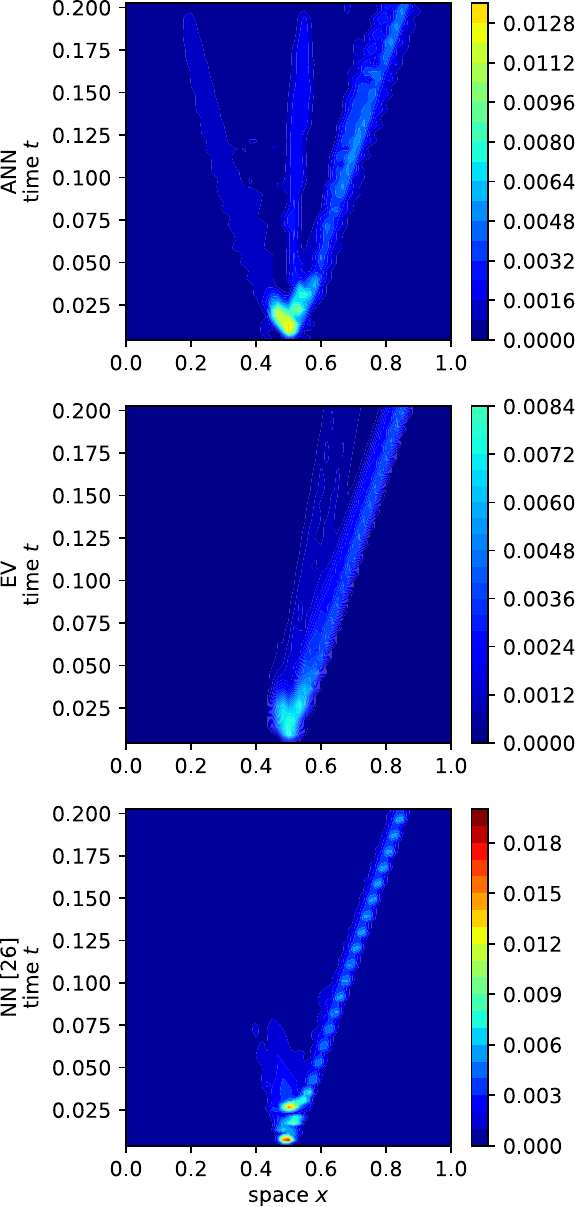}
    \end{minipage}%
    \begin{minipage}{0.32\textwidth}
        \centering
        \includegraphics[height=11.25cm]{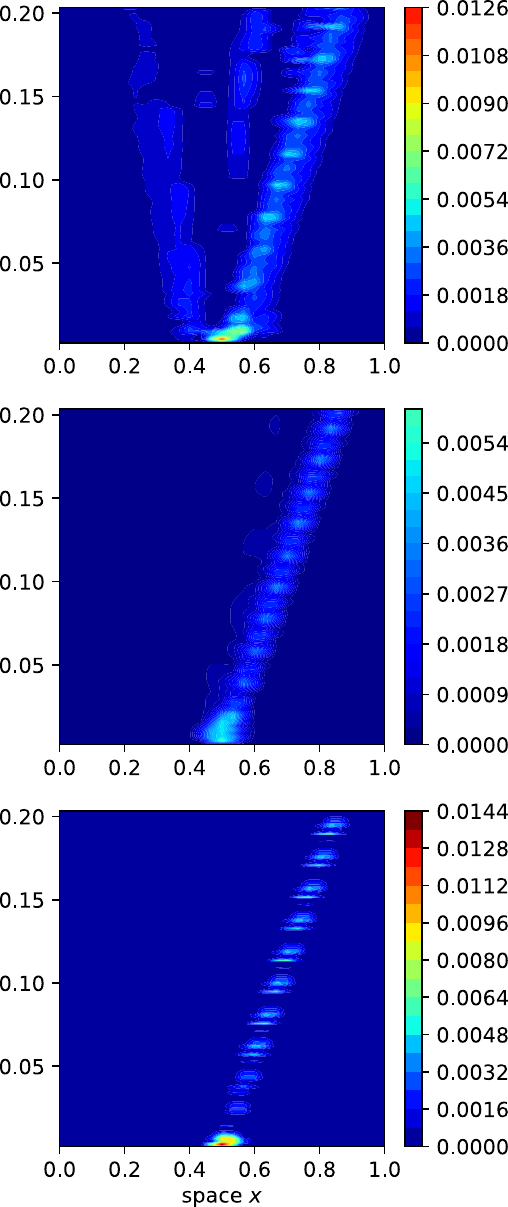}
    \end{minipage}%
    \begin{minipage}{0.32\textwidth}
        \centering
        \includegraphics[height=7.5cm]{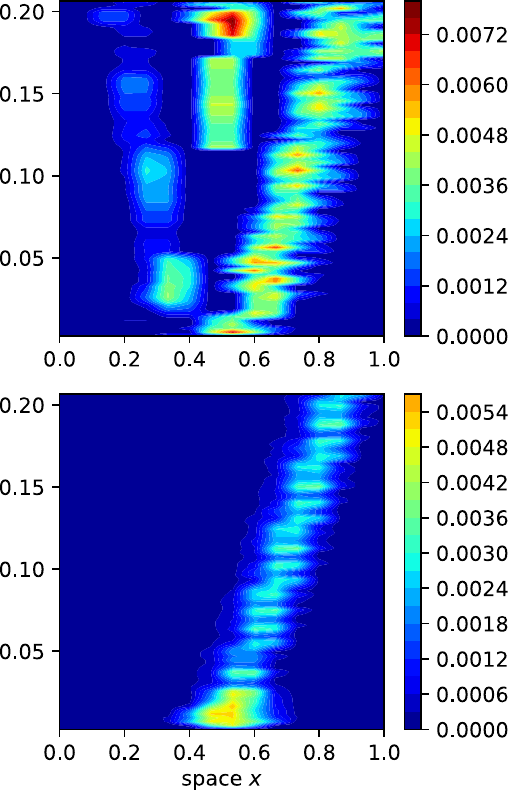}
        \vspace{3.75cm}
    \end{minipage}
    \caption{Test case 5: comparison of space-time artificial viscosity for $k=1,3,5$ (from left to right).}
    \label{fig:t5-mu}
\end{figure}


\begin{table}[t]
    \centering
    \footnotesize
    \begin{tabular}{l|ccc|ccc|cc}
                        & \multicolumn{3}{c}{$k=1$} & \multicolumn{3}{c}{$k=3$} & \multicolumn{2}{c}{$k=5$} \\
                         & NN Ours       & EV        & \textcolor{black}{NN [26]}   & NN Ours       & EV        & \textcolor{black}{NN [26]}   & NN Ours         & EV        \\
        \hline
        $\epsilon$             & 5.7046e1  & 5.8845e1  & \textcolor{black}{5.9923e1 } & 5.8787e1  & 7.2200e1  & \textcolor{black}{6.5057e1}  & 4.3138e1  & 5.7406e1 \\
        $\nabla\epsilon$       & 2.7103e1  & 2.7001e1  & \textcolor{black}{2.7639e1 } & 1.3087e2  & 1.3356e2  & \textcolor{black}{1.3892e2}  & 2.0463e2  & 2.1445e2 \\
        $\jmp{\epsilon}$       & 1.5025e0  & 1.8547e0  & \textcolor{black}{1.8261e0 } & 1.0929e0  & 1.4114e0  & \textcolor{black}{2.3391e0}  & 8.7019e-1 & 1.0149e0 \\
        o/u                    & 3.0200e-1 & 7.5084e-1 & \textcolor{black}{6.3882e-1} & 2.0631e-1 & 7.5516e-1 & \textcolor{black}{1.7252e0}  & 8.7122e-1 & 1.1302e0 \\
        \textcolor{black}{mv} &  \textcolor{black}{1.0951e-8} &  \textcolor{black}{5.0988e-8} &  \textcolor{black}{1.7928e-8} &  \textcolor{black}{2.3772e-7} &  \textcolor{black}{3.2462e-6} &  \textcolor{black}{2.9315e-6} &  \textcolor{black}{1.0951e-8} &  \textcolor{black}{5.0988e-8} \\
    \end{tabular}
    \caption{Test case 5: comparison of cumulative $L^1$ error metrics Eq.~\bref{eq:err-metric}~-~\bref{eq:ou-metric} over all timesteps from $0$ to $T$.}
    \label{tab:t5-err}
\end{table}

We now consider the first of two test cases for the Euler equations, namely the Sod shock-tube problem \cite{sod1978survey}. The problem is defined by Eq.~\bref{eq:euler-1d}, it prescribes $\Omega = (0, 1)$, $T=0.2$ and the initial state
\begin{equation*}
    \boldsymbol u_0(x) =
    \begin{bmatrix}
        \rho_0\\ v_0\\ p_0
    \end{bmatrix}
    = 
    \begin{cases}
        [1, 0, 1]^\top & x \in (0, \frac{1}{2}],\\
        [\frac{1}{8}, 0, \frac{1}{10}]^\top & x \in (\frac{1}{2}, 1).
    \end{cases}
\end{equation*}
It is endowed with constant Dirichlet boundary conditions (cf.\ \cite{sod1978survey}). For the discretization we select three cases with $k=1,3,5$, $h=\frac{1}{60},\frac{1}{30},\frac{1}{15}$, $\text{CFL}=0.27,0.61,0.88$, respectively. \textcolor{black}{Even if this problem is present in the training dataset, the polynomial degree $k$ or CFL used to discretize it are different.} \textcolor{black}{The Courant number is $C=0.2970,0.1393,0.0785$.} This is done to compare solutions with a similar number of degrees of freedom. The discontinuities in the initial conditions give rise to three characteristic waves, demanding robust numerical schemes to accurately compute the solution in the presence of shock waves and rarefactions. Namely, these are a right-moving contact wave and shock wave, and a left-moving rarefaction wave. Classical models aim to add a small amount of viscosity near the contact while constantly introducing dissipation in the region close to the shock. In Figure~\ref{fig:t5-ufin} we compare the NN solution at final time with the EV tune model with $c_K=1.0, c_{\max{}}=0.5$ \textcolor{black}{and the NN method proposed in \cite{discacciati2020controlling}. Our} neural viscosity model provides a solution that has much smaller overshoots (see $x=0.2, k=1$ for instance) and that captures with more precision the wavefront (see $x=0.7, k=1$ for instance). In this case, it is particularly interesting to analyze the history of the viscosity shown in Figure~\ref{fig:t5-mu}. Namely, here we can notice a large discrepancy between what the NN and the EV model aim to do: the EV model stabilizes only the shock, meanwhile the NN \textcolor{black}{identifies and tries} to stabilize all three present waves. This is counterbalanced by injecting a smaller dissipation near the shock. 
\textcolor{black}{For $k=5$, our model injects a large amount of viscosity near the contact discontinuity, which does not make physical sense. However, 
it is hard to identify the precise reasons why the model introduces such a large amount of artificial viscosity.
Indeed, one of the limitations of our \textit{hybrid} approach is that the model learns the viscosity from the data without explicit guidance on what constitutes correct or incorrect behavior in specific contexts.} 
Despite being more intrusive, our NN outperforms the EV model \textcolor{black}{and the NN of \cite{discacciati2020controlling}} under all the considered metrics, as shown in Table~\ref{tab:t5-err}, where we report the error metrics in time and their cumulative value.

\subsubsection{Test case 6: Shu-Osher problem}


\begin{figure}[!t]
    \centering
    \begin{minipage}{0.33\textwidth}
        \centering
        \includegraphics[width=0.95\linewidth]{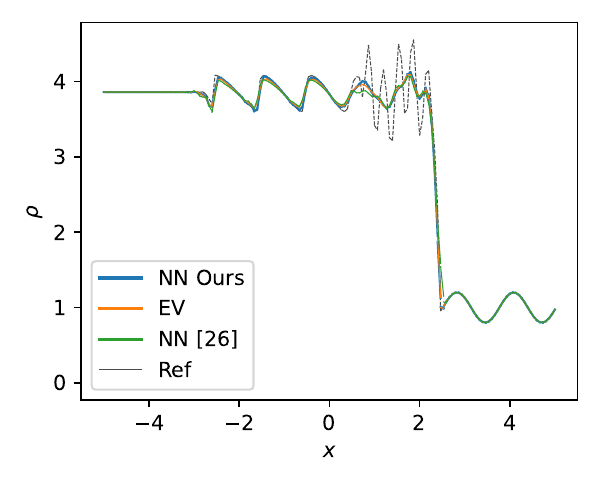}
    \end{minipage}%
    \begin{minipage}{0.33\textwidth}
        \centering
        \includegraphics[width=0.95\linewidth]{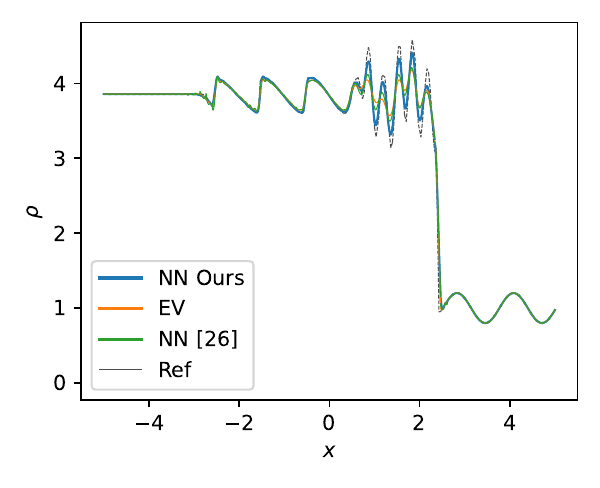}
    \end{minipage}%
    \begin{minipage}{0.33\textwidth}
        \centering
        \includegraphics[width=0.95\linewidth]{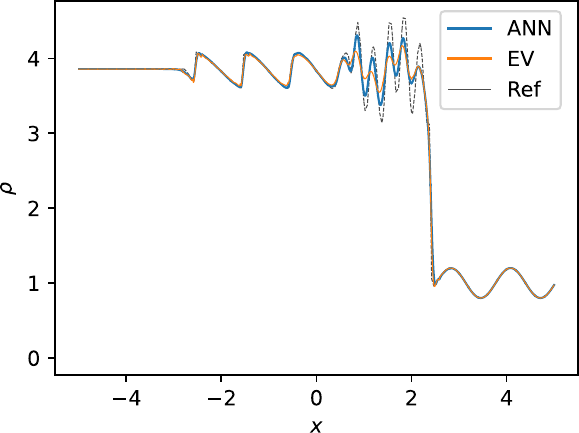}
    \end{minipage}
    \caption{Test case 6: density $\rho$ at $T=1.8$ for $k=1,3,5$ (from left to right).}
    \label{fig:t6-ufin}
\end{figure}

\begin{figure}[!t]
    \centering
    \begin{minipage}{0.36\textwidth}
        \centering
        \includegraphics[height=11.25cm]{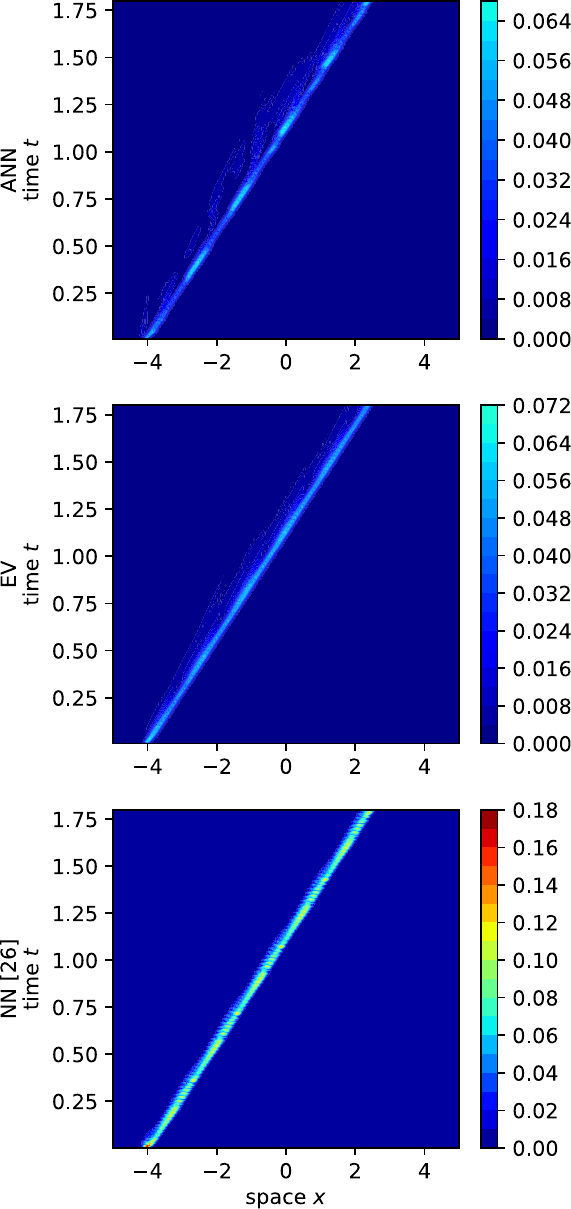}
    \end{minipage}%
    \begin{minipage}{0.32\textwidth}
        \centering
        \includegraphics[height=11.25cm]{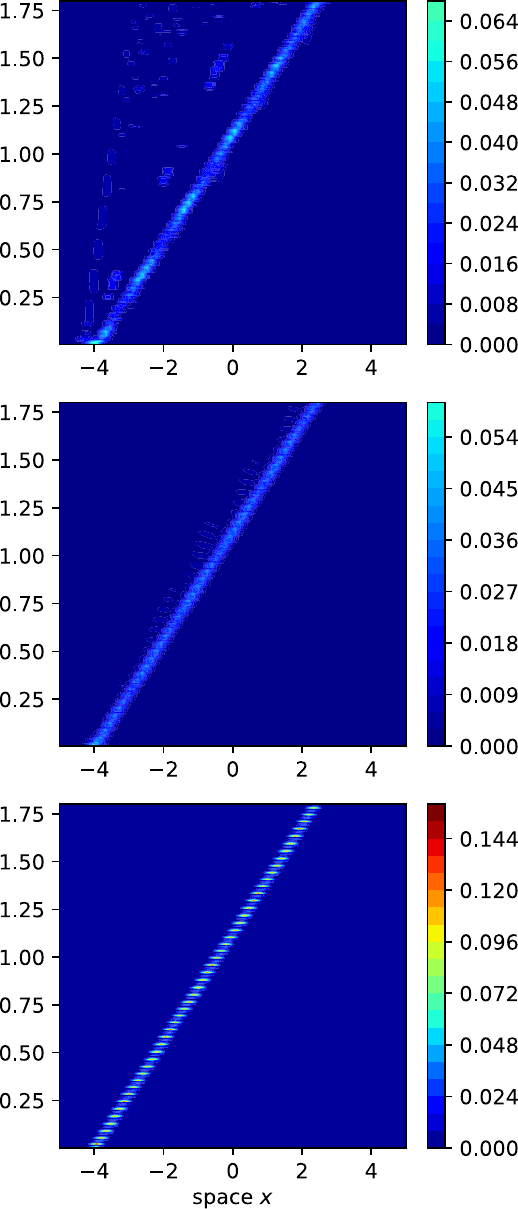}
    \end{minipage}%
    \begin{minipage}{0.32\textwidth}
        \centering
        \includegraphics[height=7.5cm]{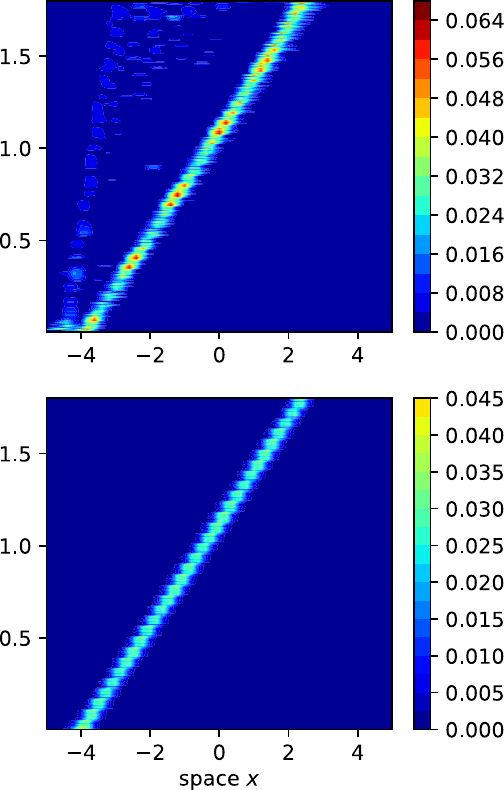}
        \vspace{3.75cm}
    \end{minipage}
    \caption{Test case 6: comparison of space-time artificial viscosity for $k=1,3,5$ (from left to right).}
    \label{fig:t6-mu}
\end{figure}
 
\begin{table}[!t]
    \centering
    \footnotesize
    \begin{tabular}{l|ccc|ccc|cc}
                        & \multicolumn{3}{c}{$k=1$} & \multicolumn{3}{c}{$k=3$} & \multicolumn{2}{c}{$k=5$} \\
                        & NN Ours      & EV       & \textcolor{black}{NN [26]}  & NN Ours       & EV      & \textcolor{black}{NN [26]}  & NN Ours       & EV        \\
        \hline
        $\epsilon$             & 2.6293e3 & 3.0574e3 & \textcolor{black}{3.6781e3} & 2.2953e3 & 4.0510e3 & \textcolor{black}{3.1186e3} & 3.7385e3 & 5.5862e3 \\
        $\nabla\epsilon$       & 1.7160e3 & 1.8731e3 & \textcolor{black}{2.0671e3} & 5.6320e3 & 7.8568e3 & \textcolor{black}{6.1255e3} & 1.3576e4 & 1.6539e4 \\
        $\jmp{\epsilon}$       & 1.2497e2 & 1.0530e2 & \textcolor{black}{8.9058e1} & 9.8602e1 & 7.1604e1 & \textcolor{black}{7.2939e1} & 5.2003e1 & 5.4259e1 \\
        o/u                    & 3.1512e0 & 3.9106e0 & \textcolor{black}{1.0716e0} & 2.8027e0 & 5.7163e0 & \textcolor{black}{2.5512e0} & 3.8023e0 & 4.6424e0 \\
        \textcolor{black}{mv} & \textcolor{black}{4.2955e-7} & \textcolor{black}{2.2112e-6} & \textcolor{black}{1.0787e-6} & \textcolor{black}{8.0499e-5} & \textcolor{black}{2.2317e-4} & \textcolor{black}{1.1017e-4} & \textcolor{black}{7.3125e-7} & \textcolor{black}{4.5271e-6} \\
    \end{tabular}
    \caption{Test case 6: comparison of cumulative $L^1$ error metrics Eq.~\bref{eq:err-metric}~-~\bref{eq:ou-metric} over all timesteps from $0$ to $T$.}
    \label{tab:t6-err}
\end{table}

The second Euler problem we consider is the Shu-Osher problem \cite{shu1988efficient}. It deals with a shock front moving inside a one-dimensional inviscid flow with artificial sinusoidal density fluctuations. This test case is relevant because it benchmarks the ability of the solver to represent both shocks and physical oscillations created by a union of smooth and discontinuous initial data. Namely, we consider the flux of Eq.~\bref{eq:euler-1d}, $\Omega = (-5, 5)$, $T=1.8$ and initial condition
\begin{equation*}
    \boldsymbol u_0(x) =
    \begin{bmatrix}
        \rho_0\\ v_0\\ p_0
    \end{bmatrix}
    = 
    \begin{cases}
        [3.857143, 2.629369, 10.333333]^\top  & x \in (-5, -4],\\
        [1 + \frac{1}{2} \sin(5x), 0, 1]^\top & x \in (-4, 5).
    \end{cases}
\end{equation*}
The problem is completed with Dirichlet boundary conditions on the left and Neumann on the right part of the boundary (cfr.\ \cite{shu1988efficient}). Since the solution exhibits a high frequency wave we employ a finer discretization, namely $h=\frac{1}{150},\frac{1}{75},\frac{1}{50}$ with $k=1,3,5$, and $\text{CFL}=0.12,0.3,0.4$, respectively. \textcolor{black}{The Courant number is $C=0.2707,0.1337,0.0762$.} From Figure~\ref{fig:t6-ufin}, we can see that for the case $k=1$ the NN and the EV model \textcolor{black}{and the NN of \cite{discacciati2020controlling}} produce a rather similar solution. However, the quantitative analysis reported in Table~\ref{tab:t6-err} shows that the neural viscosity model is slightly more accurate in the $L^1$ norm, and when considering overshoots and undershoots. Instead, the cases $k=3,5$ show that the EV model \textcolor{black}{and the NN of \cite{discacciati2020controlling}} smooths the solution too much, and thus it is not able to capture the high frequencies. On the other hand, \textcolor{black}{our} NN shows much better results, namely the solution correctly represents the high-frequency wave pack present between $x=0.5$ and $x=2.5$. This also proves that the neural viscosity model is able to exploit the high-order nature of the method, indeed the three considered cases have the same number of degrees of freedom. In particular, for $k=3$ we have a reduction of the $L^1$ error of 43\% and of 51\% of overshoot/undershoot. The dissipation injected by the NN, as shown in Figure~\ref{fig:t6-mu}, is concentrated on the shock and, differently from the EV model, adds also other small patches of viscosity at the boundary of the sinusoidal wave pack.

\subsection{Two-dimensional problems}
We analyze now two-dimensional problems. Notice that, to pass to 2D, we do not have to change anything about the training algorithm or the neural network architecture: everything is exactly the same as in the one-dimensional case. However, the physics that governs these problems is significantly different from before. This requires to train a new neural network, which in turn entails a re-tuning of the hyperparameters. Indeed, in two dimensions more complex structures appear as a result of the interactions among one-dimensional waves. Hence, the differences among the models are more prominent, requiring an NN able to generalize better.

We choose the hyperparameter with a procedure that exactly mimics what was done in the one-dimensional case, as described in Subsection~\ref{sec:training-in-practice}. Let us remark that this expensive operation of training and hyperparameters tuning must be done just once for all the two-dimensional problems, meanwhile, traditional models must be tuned problem by problem. At the end of this process, we selected a NN with a width of 160 and a depth of eight (the other hyperparameters remain unchanged). As we have done for the one-dimensional test cases, we start by showcasing the results on scalar problems and then pass to the Euler system of equations.

\subsubsection{Test case 7: linear advection}

\begin{figure}[t]
    \centering
    \includegraphics[width=0.5\linewidth]{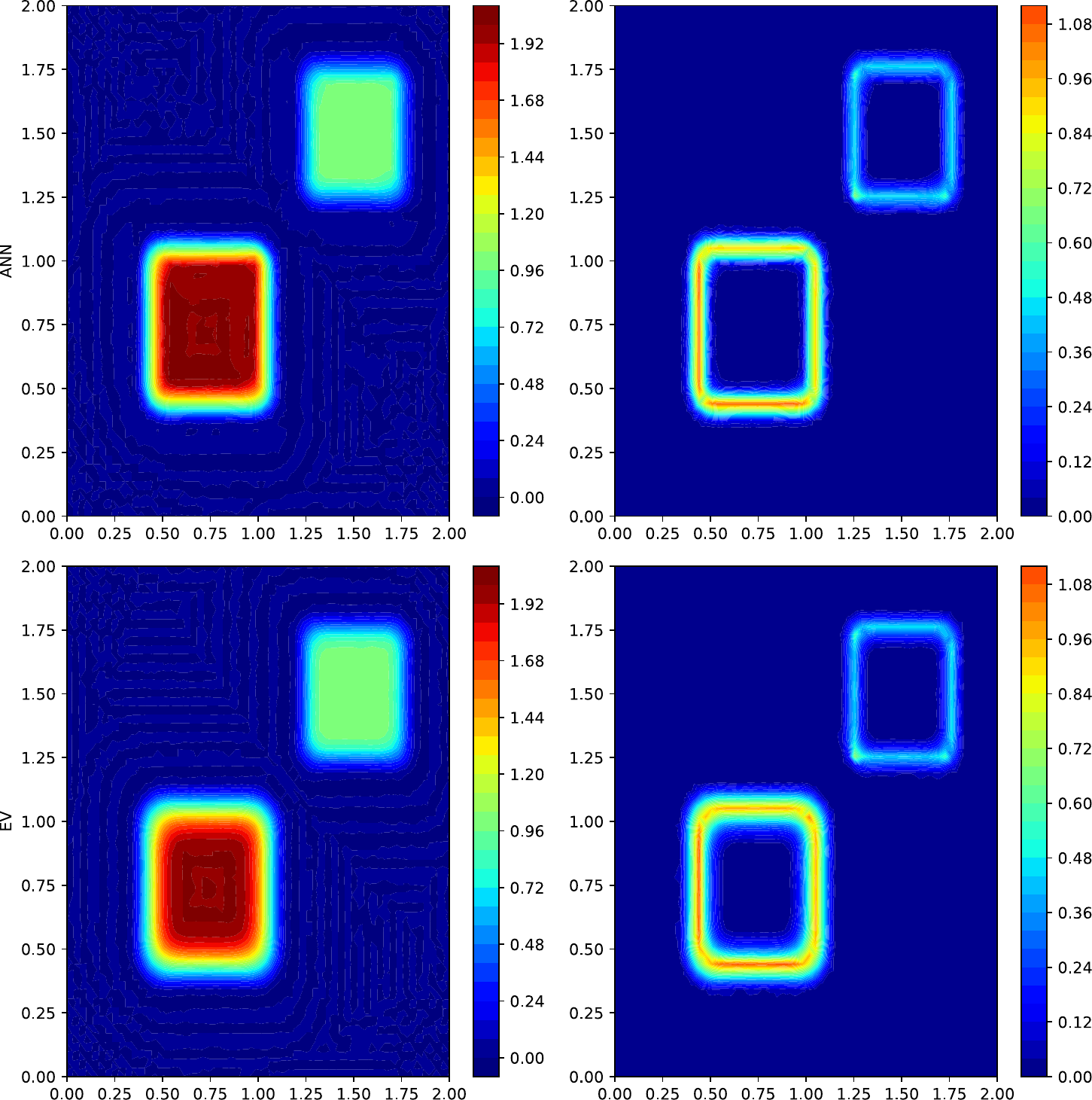}
    \caption{Test case 7: solution and error at final time $T=0.25$ for $k=1$.}
    \label{fig:t7-ufin}
\end{figure}

\begin{figure}[t]
    \centering
        \includegraphics[width=0.78\textwidth,trim={0 0 0 0},clip]{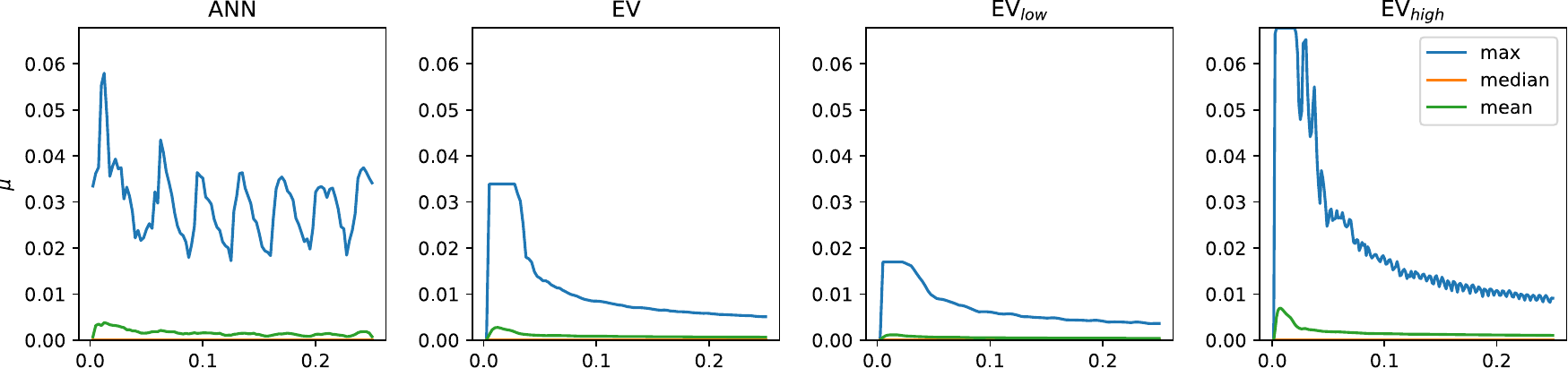}
    \caption{Test case 7: maximum value of the artificial viscosity $\mu$ in time for $k=1$.}
    \label{fig:t7-err+mu}
\end{figure}

\begin{table}[t]
    \centering
    \footnotesize
    \begin{tabular}{c|cc|cc|cc}
                        & \multicolumn{2}{c}{$k=1$} & \multicolumn{2}{c}{$k=3$} & \multicolumn{2}{c}{$k=5$} \\
                        & NN         & EV        & NN         & EV        & NN         & EV        \\
        \hline
        $\epsilon$       & 1.0374e+05 & 1.2616e+05 & 4.7333e+04 & 6.6665e+04 & 7.3585e+04 & 8.7204e+04 \\
        $\nabla\epsilon$ & 1.0383e+05 & 1.0916e+05 & 3.0168e+05 & 3.4178e+05 & 8.6089e+05 & 8.7469e+05 \\
        $\jmp{\epsilon}$ & 9.5542e+03 & 1.0781e+04 & 4.0133e+03 & 4.7451e+03 & 2.4605e+03 & 2.3712e+03 \\
        o/u              & 8.9123e+02 & 1.6955e+03 & 2.0373e+03 & 1.1930e+03 & 5.0116e+02 & 6.5465e+02 \\
        mv               & 1.2287e-15 & 1.3516e-14 & 8.9032e-14 & 9.8162e-13 & 2.3507e-13 & 2.5857e-13 \\
    \end{tabular}
    \caption{Test case 7: comparison of cumulative $L^1$ error metrics Eq.~\bref{eq:err-metric}~-~\bref{eq:mv-metric} over all timesteps from $0$ to $T$.}
    \label{tab:t7-err}
\end{table}

The first 2D test case we consider is linear advection $\boldsymbol f(u) = \boldsymbol{\beta} u$ with coefficient $\boldsymbol{\beta} = [1, 1]^\top$. This simple test case is non-trivial when dealing with discontinuous initial condition, namely, we employ
\begin{equation*}
    u_0(x_1, x_2) =
    \begin{cases}
        2 & x_i \in (\frac{1}{5}, \frac{4}{5}), \: i\in\{1,2\},\\
        1 & x_i \in (1, \frac{3}{2}), \: i\in\{1,2\},\\
        0 &\textnormal{otherwise.}
    \end{cases}
\end{equation*}
Moreover, it gives us the possibility to check the mass conservation properties of the neural viscosity model. We solve the problem in $\Omega=(0, 2)^2$ with a final time $T=0.25$ on a structured grid with $7200, 1800, 450$ triangular elements and $k=1,3,5$, $\text{CFL}=0.075,0.3,0.32$, respectively. \textcolor{black}{The Courant number is $C=0.1061,0.0471,0.0181$.} We focus on $k=1$ since the other two cases are analogous. As done in the previous test cases, we compare the neural viscosity model with a tuned EV model ($c_K=1.0, c_{\max{}}=0.5$) and two other EV models that inject a sub-optimal amount of dissipation (too large $c_K=2.0, c_{\max{}}=1.0$ and too small $c_K=0.5, c_{\max{}}=0.25$). From Figure~\ref{fig:t7-ufin} we clearly see that both the NN and the EV model smooth the discontinuous profile of the initial condition in order to avoid oscillations. The NN solution is symmetric and visually exhibits a sharper wavefront. This is confirmed by the quantitative analysis of the error reported in Table~\ref{tab:t7-err}. In particular, it is interesting to observe that the NN model injects more viscosity than a classical EV model but has a smaller $L^1$ error. This is not an intuitive behavior when comes to the EV model since from both Figure~\ref{fig:t7-err+mu} and Table~\ref{tab:t7-err} it is clearly observable that EV models that inject larger dissipation have smaller overshoots/undershoots but larger $L^1$ error. Moreover, the solution obtained with the NN conserves the mass better than the one obtained with the EV model, showing that also physics is respected. Therefore, in this case, the neural viscosity model is effectively learning a new model that has better properties than the EV one.

\subsubsection{Test case 8: KPP rotating wave}

\begin{figure}[t]
    \centering
    \begin{minipage}{0.5\textwidth}
        \centering
        \includegraphics[width=0.85\linewidth]{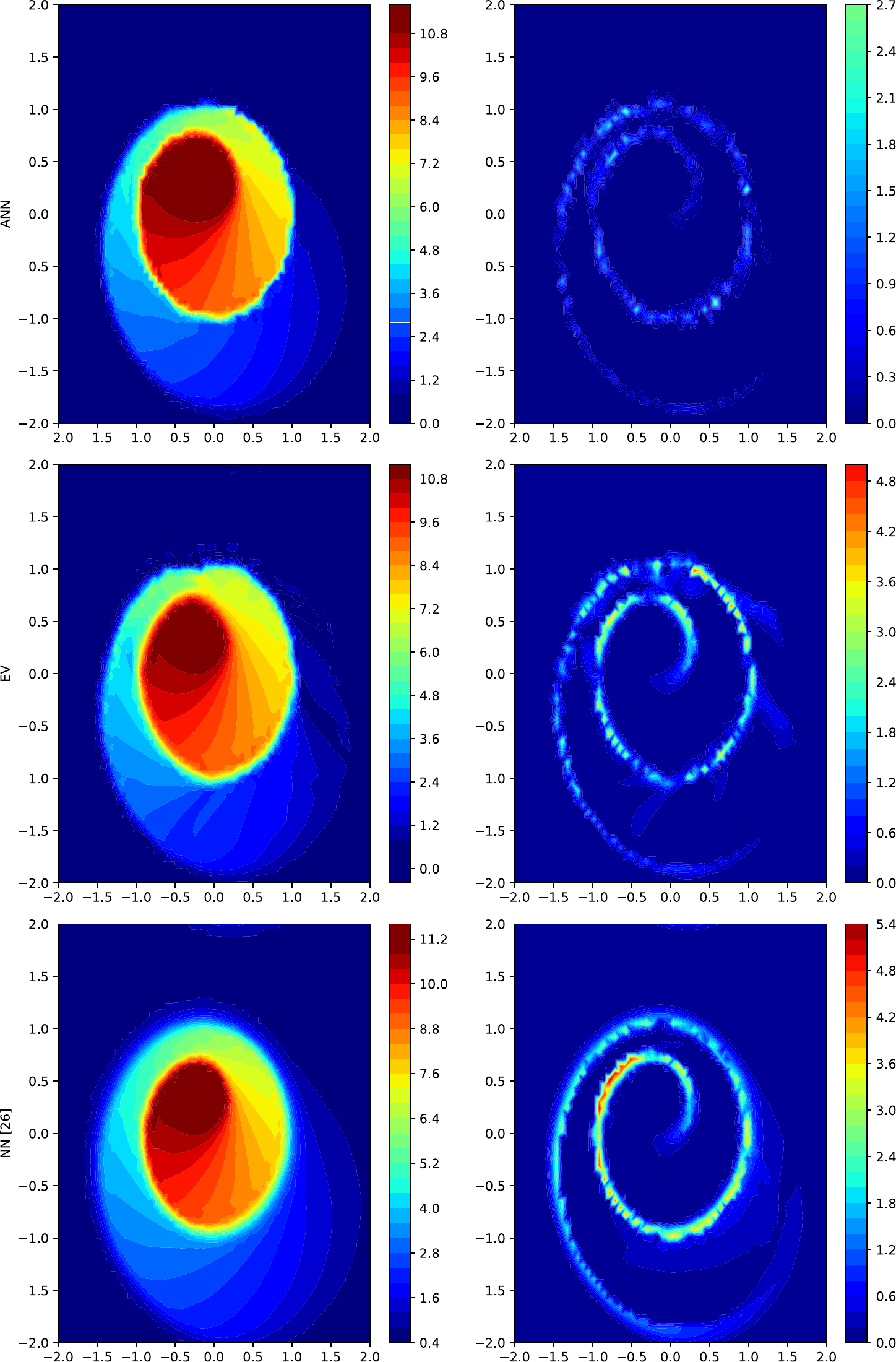}
    \end{minipage}%
    \begin{minipage}{0.5\textwidth}
        \centering
        \includegraphics[width=0.85\linewidth]{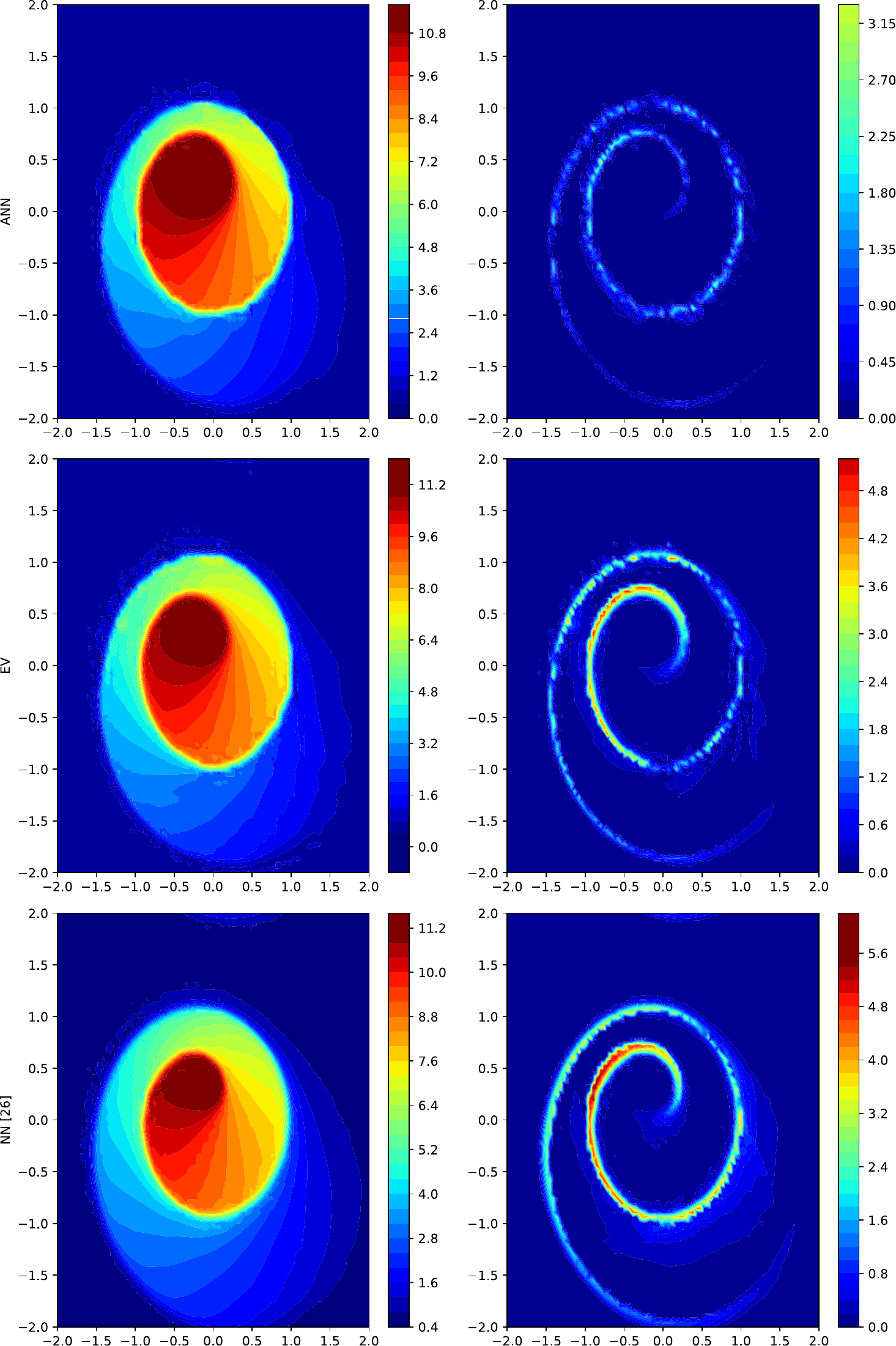}
    \end{minipage}%
    \caption{Test case 8: solution and error at final time $T=0.25$ for $k=1,3$.}
    \label{fig:t8-k13-sol}
\end{figure}

\begin{figure}[t]
        \centering
        \includegraphics[width=0.4\linewidth]{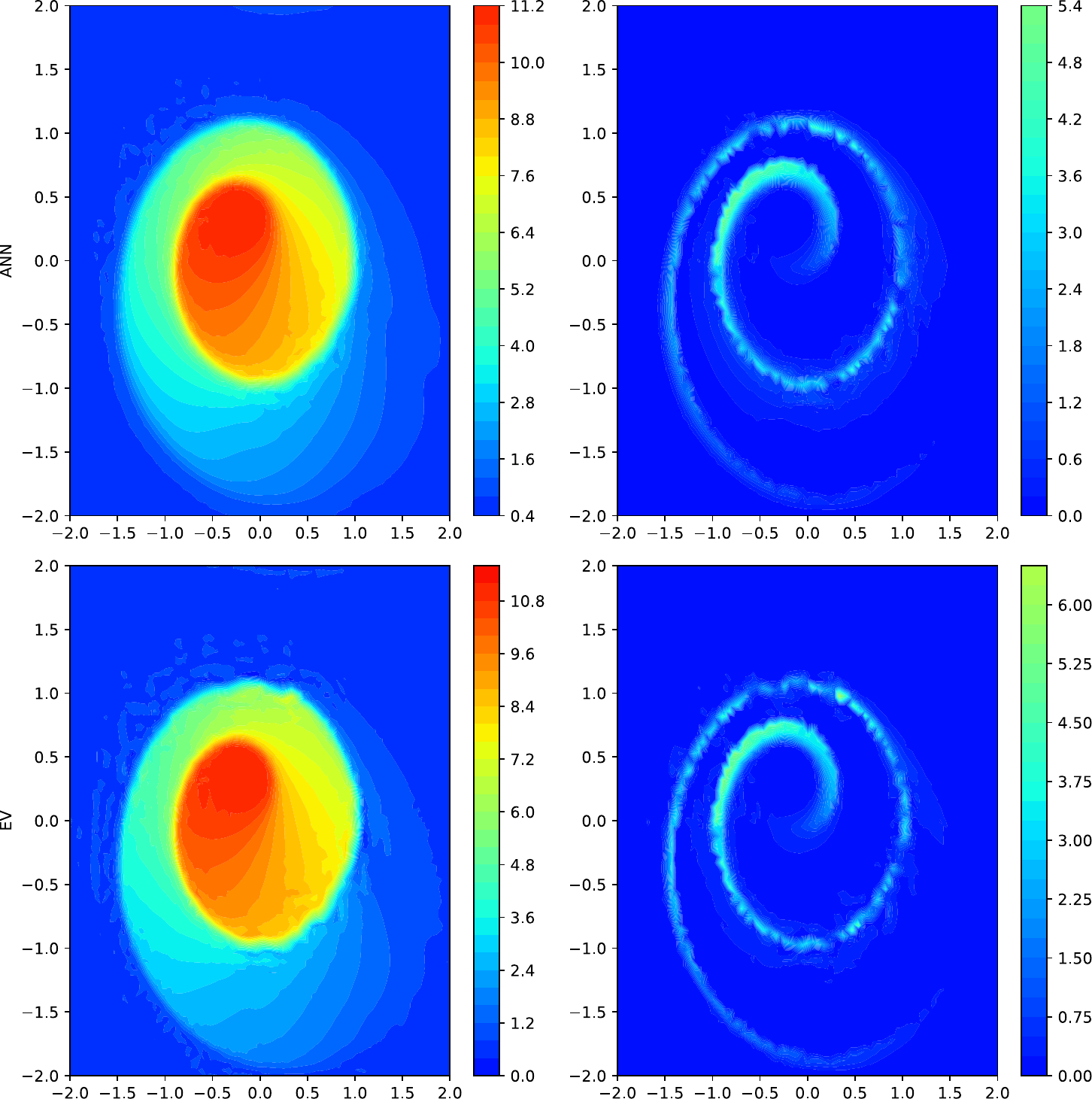}
        \caption{Test case 8: solution and error at final time $T=0.25$ for $k=5$.}
        \label{fig:t8-k5-sol}
\end{figure}

\begin{figure}[t]
    \centering
    \includegraphics[width=0.85\textwidth,trim={0 0.9cm 0 0},clip]{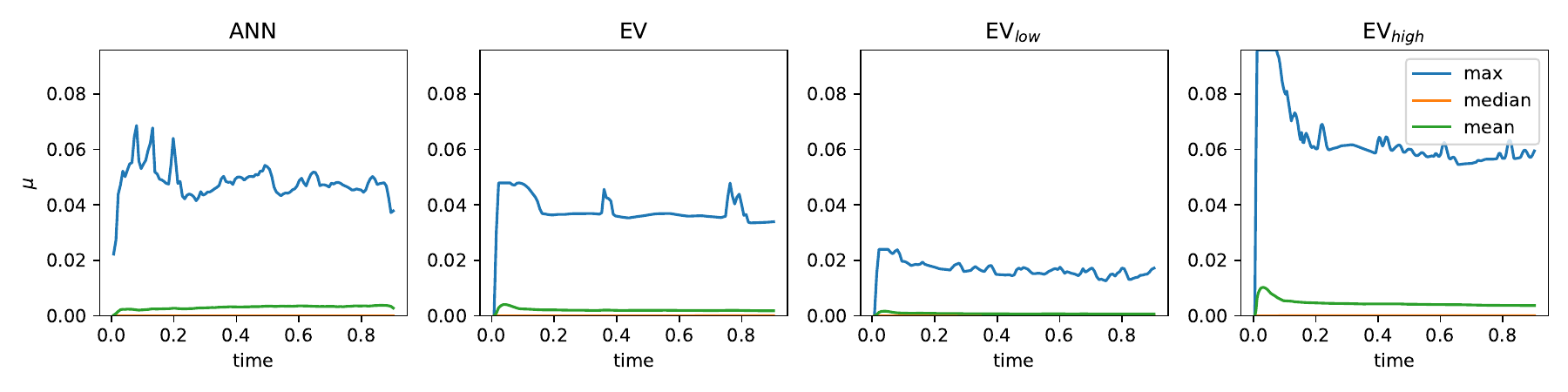}
    \includegraphics[width=0.85\textwidth,trim={0 0.9cm 0 0.9cm},clip]{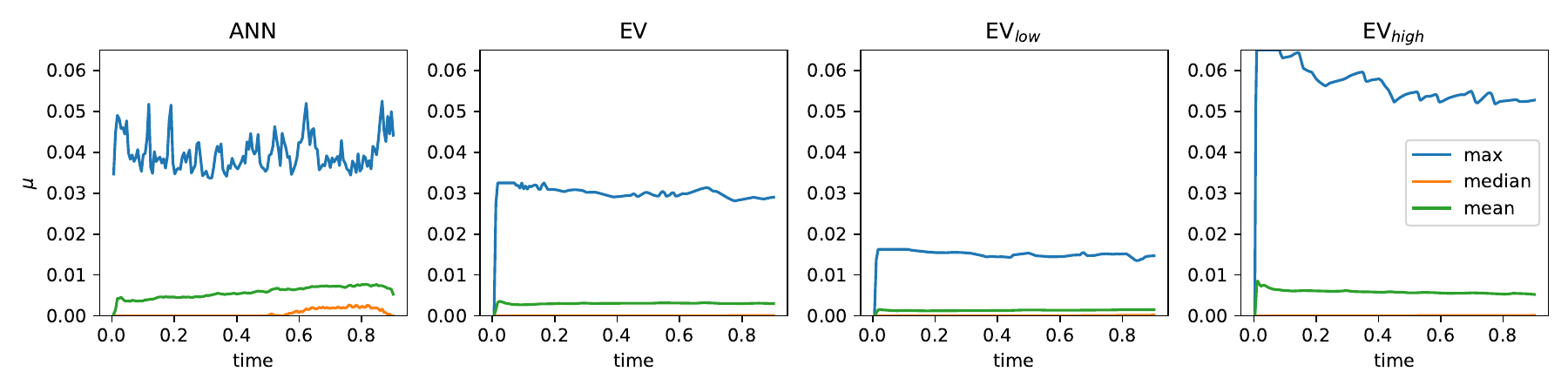}
    \includegraphics[width=0.85\textwidth,trim={0 0 0 0.9cm},clip]{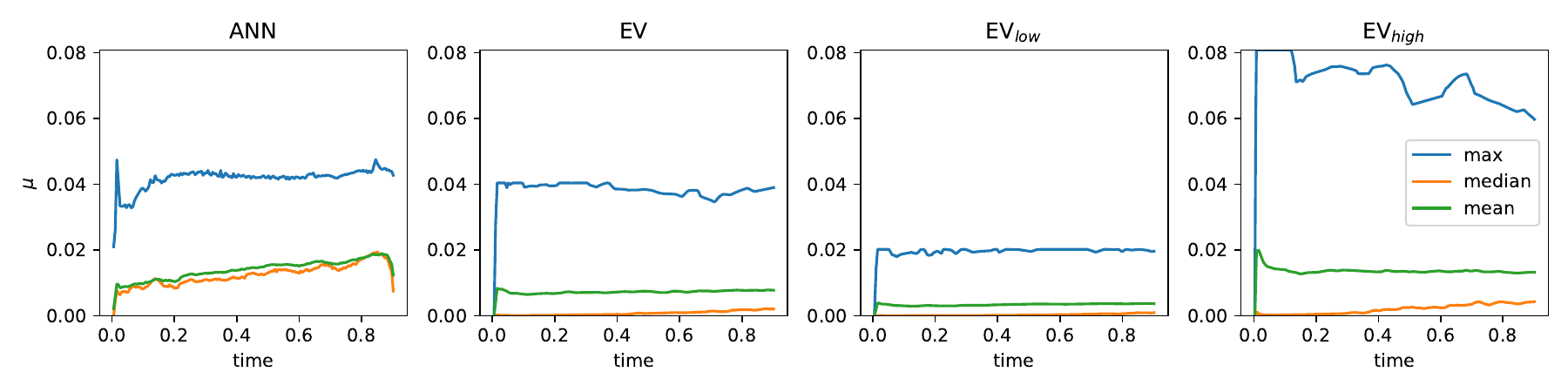}
    \caption{Test case 8: maximum value of the artificial viscosity $\mu$ in time for $k=1,3,5$ (from top to bottom).}
    \label{fig:t8-mu}
\end{figure}

\begin{table}[t]
    \centering
    \footnotesize
    \begin{tabular}{c|ccc|ccc|cc}
                        & \multicolumn{3}{c}{$k=1$} & \multicolumn{3}{c}{$k=3$} & \multicolumn{2}{c}{$k=5$} \\
                         & NN Ours  & EV       & \textcolor{black}{NN [26]}  & NN       & EV       & \textcolor{black}{NN [26]}  & NN       & EV       \\
        \hline
        $\epsilon$       & 1.6196e5 & 3.9369e5 & \textcolor{black}{6.3515e5} & 1.6956e5 & 3.8346e5 & \textcolor{black}{5.6881e5} & 2.5583e5 & 2.8887e5 \\
        $\nabla\epsilon$ & 2.4551e5 & 4.3755e5 & \textcolor{black}{5.4866e5} & 1.2625e6 & 1.7127e6 & \textcolor{black}{1.8699e6} & 2.6152e6 & 2.8368e6 \\
        $\jmp{\epsilon}$ & 1.0973e5 & 1.4006e5 & \textcolor{black}{9.5575e4} & 8.2736e4 & 9.7502e4 & \textcolor{black}{6.4618e4} & 4.4259e4 & 6.2624e4 \\
        o/u              & 8.0091e3 & 1.9289e4 & \textcolor{black}{4.8742e3} & 6.4779e3 & 6.8667e3 & \textcolor{black}{2.5480e3} & 4.6155e3 & 5.0963e3 \\
    \end{tabular}
    \caption{Test case 8: comparison of cumulative $L^1$ error metrics Eq.~\bref{eq:err-metric}~-~\bref{eq:ou-metric} over all timesteps from $0$ to $T$.}
    \label{tab:t8-err}
\end{table}

We consider now nonlinear scalar conservation equations characterized by the non-convex flux $[\sin u, \cos u]^\top$. Specifically, we examine a two-dimensional scalar conservation equation introduced in \cite{guermond2011entropy, kurganov2007adaptive}. This particular test poses a challenge to numerous high-order numerical schemes due to its intricate two-dimensional composite wave structure in the solution. Namely, we have $\Omega=(-2, 2)^2$, $T=1$ and initial condition

\begin{equation*}
    u_0(\boldsymbol x) = 
    \begin{cases}
        \frac{7}{2}\pi & \textnormal{ if } \norm{\boldsymbol x}_2 < 1,\\
        \frac{1}{4}\pi & \textnormal{ otherwise.}
    \end{cases}
\end{equation*}
The problem is endowed with periodic boundary conditions. For discretization, we employ a structured triangular mesh with $7200, 1800, 450$ elements and $k=1,3,5$, $\text{CFL}=0.11,0.4,0.48$, respectively. \textcolor{black}{The Courant number is $C=0.11,0.0444,0.0192$.} We compare the neural viscosity model with EV models tuned like in the previous test case \textcolor{black}{and the NN of \cite{discacciati2020controlling}}. The results reported in Figure~\ref{fig:t8-k13-sol} and show that the neural viscosity model produces more accurate solutions with respect to the EV one. In particular, it presents less wriggles and a sharper wave front. To understand the behavior of the model we report also in Figure~\ref{fig:t8-mu} the maximum, mean, and median of the injected viscosity. As was the case in 1D, the NN model injects more viscosity than the tuned EV model, but retains a lower $L^1$ error, as shown in Table~\ref{tab:t8-err}. \textcolor{black}{On the other hand, the method proposed in \cite{discacciati2020controlling}, is the one that produces the smoothest solution, featuring far less oscillations than the other two methods. However, it is also the less accurate solution from the point of view of the $L^1$ norm.}

\subsubsection{Test case 9: Euler equation, Riemann problem}
 
\begin{figure}[!t]
    \begin{minipage}{0.5\textwidth}
        \centering
        \includegraphics[width=0.95\linewidth]{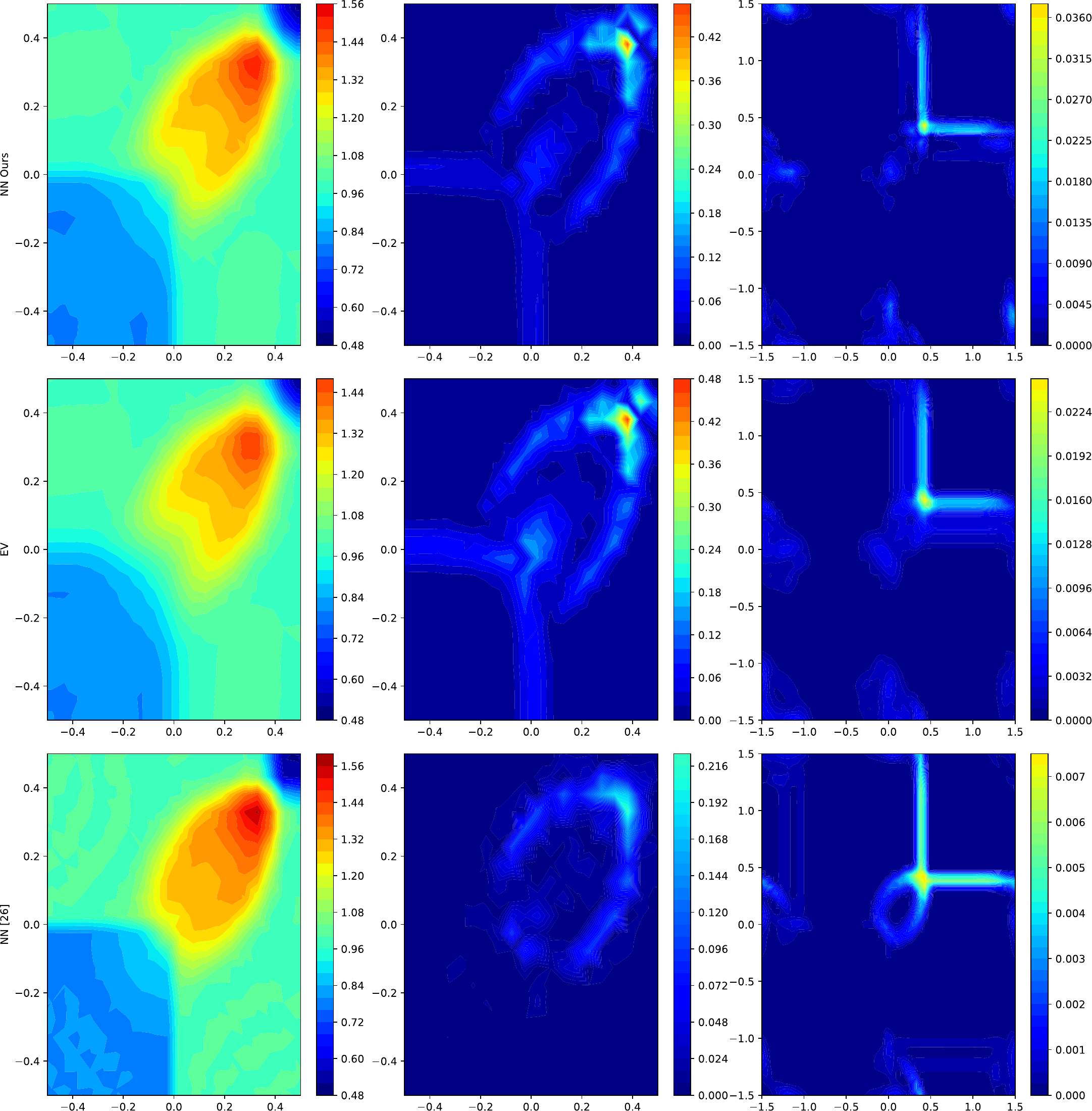}
    \end{minipage}%
    \begin{minipage}{0.5\textwidth}
        \centering
        \includegraphics[width=0.95\linewidth]{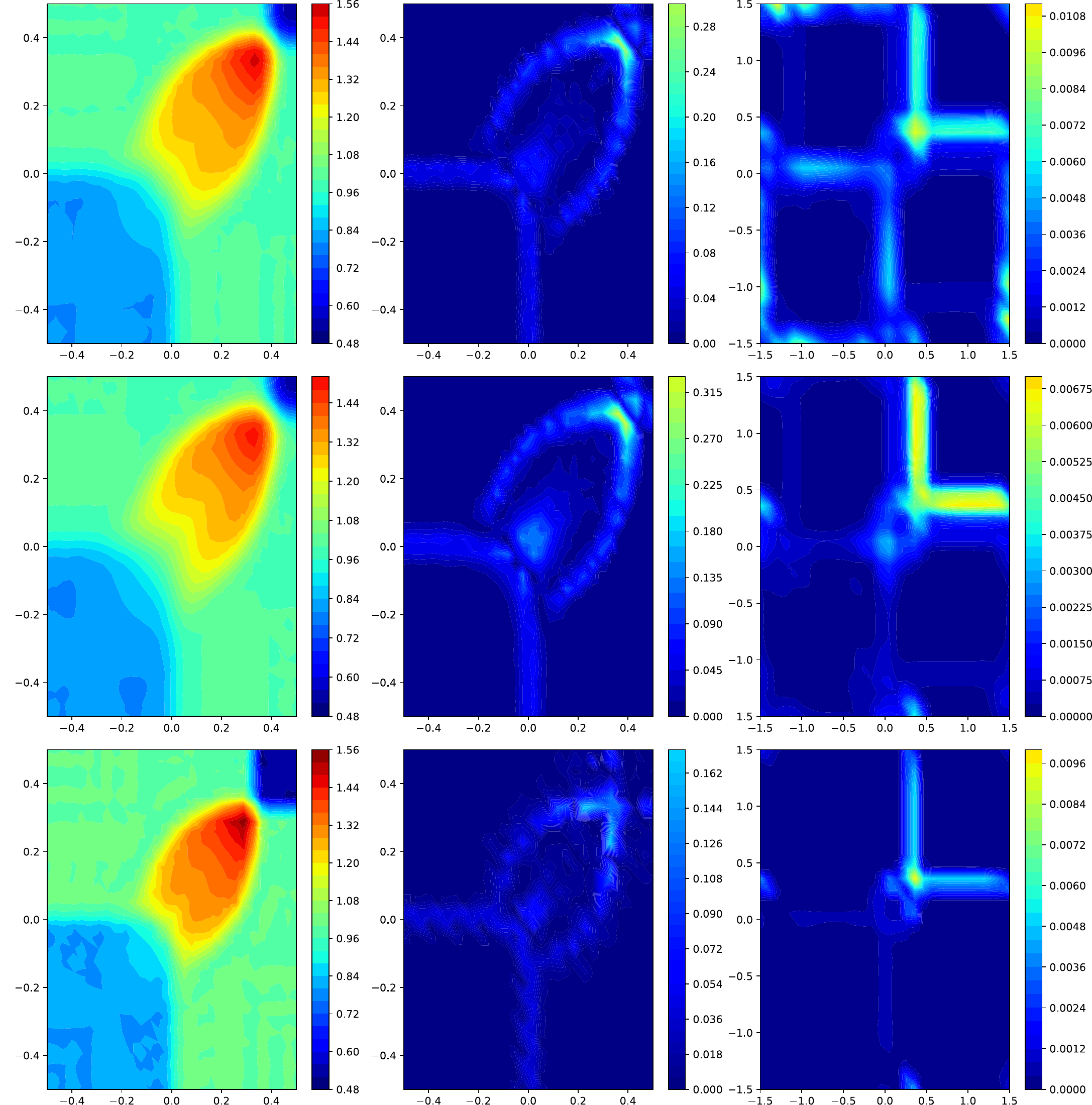}
    \end{minipage}%
    \caption{Test case 9: solution and error at final time $T=0.25$ for $k=1,3$ (from left to right).}
    \label{fig:t9-ufin-k13}
\end{figure}

\begin{figure}[!t]
        \centering
        \includegraphics[width=0.4\linewidth]{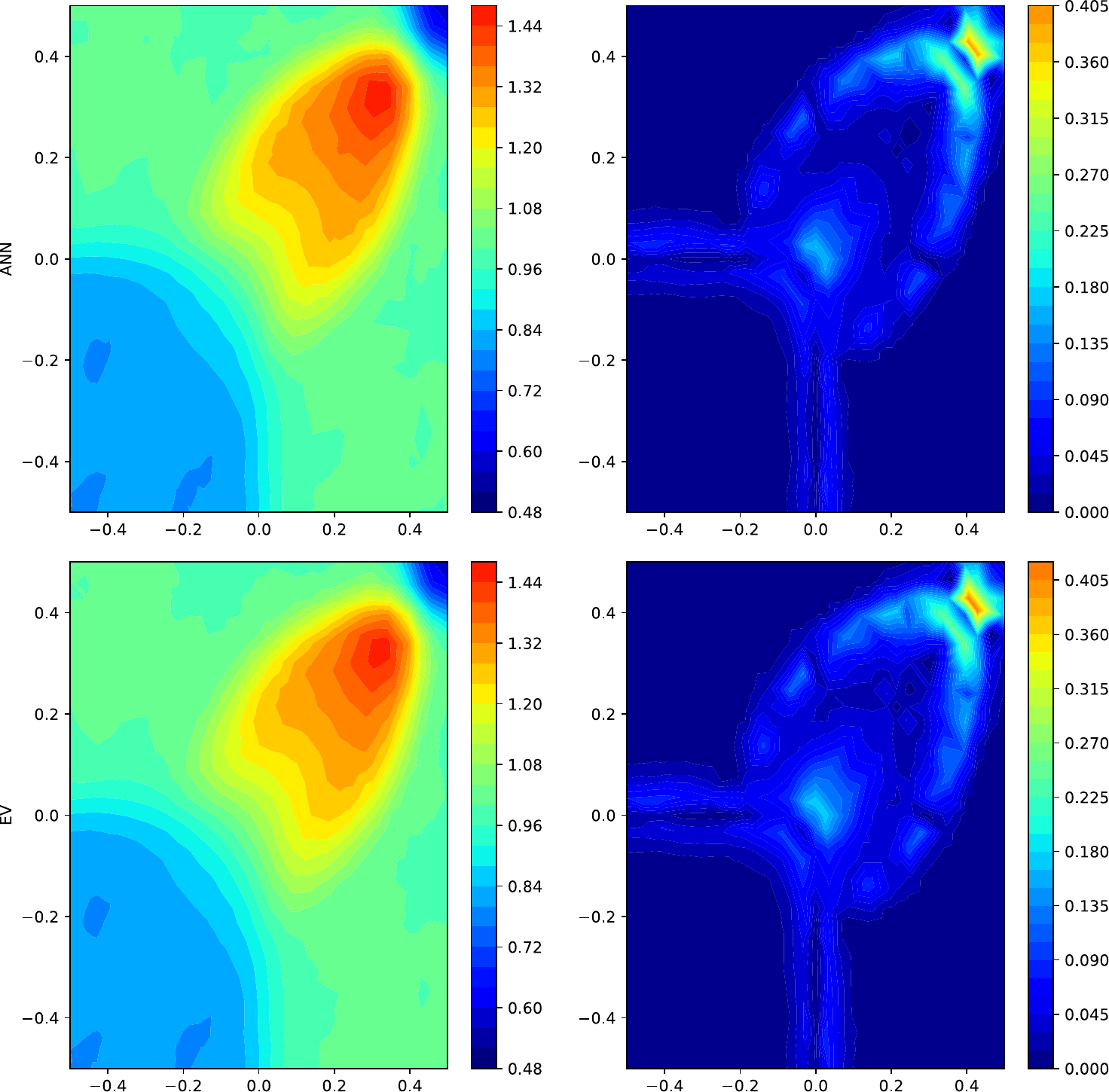}
        \caption{Test case 9: solution and error at final time $T=0.25$ for $k=5$.}
        \label{fig:t9-ufin-k5}
\end{figure}


\begin{table}[t]
    \centering
    \footnotesize
    \begin{tabular}{c|ccc|ccc|cc}
                        & \multicolumn{3}{c}{$k=1$} & \multicolumn{3}{c}{$k=3$} & \multicolumn{2}{c}{$k=5$} \\
                         & NN Ours  & EV       & \textcolor{black}{NN [26]}  & NN Ours  & EV       & \textcolor{black}{NN [26]}  & NN Ours         & EV        \\
        \hline
        $\epsilon$             & 1.8409e4 & 3.1180e4 & \textcolor{black}{1.4657e4} & 1.2387e4 & 2.1667e4 & \textcolor{black}{1.1576e4} & 1.1578e4 & 1.7845e4 \\
        $\nabla\epsilon$       & 2.0193e4 & 2.9233e4 & \textcolor{black}{1.6856e4} & 7.0429e4 & 9.1563e4 & \textcolor{black}{6.6923e4} & 1.3498e5 & 1.6098e5 \\
        $\jmp{\epsilon}$       & 5.5618e3 & 4.5949e3 & \textcolor{black}{7.2912e3} & 2.3759e3 & 1.7313e3 & \textcolor{black}{4.4055e3} & 1.4495e3 & 9.8379e2 \\
        o/u                    & 7.4932e1 & 1.7049e2 & \textcolor{black}{6.5821e2} & 6.5138e1 & 1.3598e2 & \textcolor{black}{5.0583e2} & 3.9249e1 & 7.5978e1 \\
        \textcolor{black}{mv} & \textcolor{black}{5.911e-14} & \textcolor{black}{5.570e-14} & \textcolor{black}{5.656e-14} & \textcolor{black}{6.366e-14} & \textcolor{black}{6.480e-14} & \textcolor{black}{4.768e-14} & \textcolor{black}{7.731e-14} & \textcolor{black}{7.162e-14} \\
    \end{tabular}
    \caption{Test case 9: comparison of cumulative $L^1$ error metrics Eq.~\bref{eq:err-metric}~-~\bref{eq:ou-metric} over all timesteps from $0$ to $T$.}
    \label{tab:t9-err}
\end{table}

We now consider a two-dimensional Riemann problem for the Euler system. Namely, we consider the configuration 12 proposed in \cite{kurganov2002solution, schulz1993classification}. The problem is particularly challenging in view of the presence of both contact waves and shocks. The physical domain is $\Omega = (-\frac{1}{2}, \frac{1}{2})^2$, the final time is $T=0.25$ and the initial conditions are:
\begin{equation*}
    \boldsymbol u_0(x_1, x_2) = 
    \begin{bmatrix}
        \rho_0\\
        v_{1,0}\\
        v_{2,0}\\
        e_0
    \end{bmatrix} =
    \begin{cases}
        [0.5313, 0, 0, 0.4]^\top & \textnormal{ if } x_1 > 0, x_2 > 0,\\
        [1, 0.7276, 0, 1]^\top   & \textnormal{ if } x_1 < 0, x_2 > 0,\\
        [0.8, 0, 0, 1]^\top      & \textnormal{ if } x_1 < 0, x_2 < 0,\\
        [1, 0, 0.7276, 1]^\top   & \textnormal{ if } x_1 > 0, x_2 < 0.
    \end{cases}
\end{equation*}

Since the problem is completed with periodic boundary conditions, we enlarge the domain to $\Omega = (-\frac{3}{2}, \frac{3}{2})^2$ to avoid contaminations of the solution. The discretization is done with a structured grid with $7200, 1800, 450$ triangular elements and $k=1,3,5$, $\text{CFL}=0.05,0.18,0.21$, respectively. \textcolor{black}{The Courant number is $C=0.1098,0.0443,0.0179$.} \textcolor{black}{We compare our NN with a tuned EV model and the NN of \cite{discacciati2020controlling}}. The results are reported in Figure~\ref{fig:t9-ufin-k13}. A qualitative way to evaluate the model is to check the error done in representing the fine structures in $\boldsymbol x = [0, 0]^\top$ and $\boldsymbol x = [0.4, 0.4]^\top$. It is possible to notice that the NN model captures better the contact wave, especially for $k=3$. This is confirmed by a quantitative analysis of the errors, as shown in Table~\ref{tab:t9-err}. \textcolor{black}{In this case, the model of \cite{discacciati2020controlling} slightly outperforms our NN, producing more accurate results under several metrics.} \textcolor{black}{Our NN model injects viscosity greater than $10^{-8}$ in 3730, 1023, and 316 cells for $k=1,3,5,$ respectively.}

\section{Conclusions}\label{sec:conclusions}

In this work, we proposed a novel algorithm to train a NN surrogating an artificial viscosity model in a non-supervised way. Seeking to surpass the performance of classical viscosity models, we introduced a new approach inspired by RL that, by interacting with a DG solver, is able to discover new viscosity models without the necessity of explicitly building a dataset.

First, we presented a NN model with a renewed architecture with respect to the one present in literature \cite{discacciati2020controlling}. The modification we have introduced proved to increase the flexibility and accuracy of the network. Then, we tested the model trained with the proposed training algorithm on several test cases, ranging from problems with smooth to discontinuous solutions and from scalar to vectorial equations. The trained NN exhibits favorable generalization properties despite being trained in an environment consisting of relatively simple problems. The resulting model proves to maintain the expected order of convergence for smooth problems. Furthermore, it outperforms the optimally tuned classical models under the considered metrics. Indeed, we have demonstrated that a viscosity model with consolidated reliability, such as EV, may falter in generating accurate results unless optimal parameters are carefully selected. The trained model exhibits effective shock-capturing properties and identifies areas requiring dissipation of numerical oscillations. This feature is more apparent in two-dimensional problems, where the NN can accurately capture multidimensional waves, fine patterns, and intricate spatial configurations.

This algorithm also opens the possibility of training models by incorporating data that comes from the physical measurement of quantities of interest, such as the velocity or the pressure of a flow field. In future works, we aim to test that this integration effectively enhances model reliability by incorporating fundamental principles and constraints, leading to more accurate predictions and improved generalization to real-world scenarios. \textcolor{black}{For instance, integrating in the loss a regularization term based on the entropy production could potentially improve the model's adherence to physical laws.}
Moreover, we plan to apply the proposed algorithm to more domain-specific test cases such as the variational multiscale stabilization methods for the approximation of the incompressible Navier-Stokes equations.

\paragraph*{Declaration of competing interests.} 
The authors declare that they have no known competing financial interests or personal relationships that could have appeared to influence the work reported in this paper.

\paragraph*{Funding.}
M.C., P.F.A. and L.D. are members of the INdAM Research group GNCS. The authors have been partially funded by the research projects PRIN 2020 (n. 20204LN5N5) funded by Italian Ministry of University and Research (MUR). P.F.A. is partially funded by the European Union (ERC SyG, NEMESIS, project number 101115663). Views and opinions expressed are however those of the authors only and do not necessarily reflect those of the European Union or the European Research Council Executive Agency. Neither the European Union nor the granting authority can be held responsible for them. L.D. acknowledges the support by the FAIR (Future Artificial Intelligence Research) project, funded by the NextGenerationEU program within the PNRR-PE-AI scheme (M4C2, investment 1.3, line on Artificial Intelligence), Italy. The present research is part of the activities of “Dipartimento di Eccellenza 2023-2027”, MUR, Italy.

\paragraph*{CRediT authorship contribution statement.}
P.F.A. and L.D.: Conceptualization, Methodology, Resources, Review and editing, Supervision, Project administration, Funding acquisition. M.C.: Conceptualization, Methodology, Software, Validation, Formal analysis, Investigation, Data curation, Visualization, Original draft.

\bibliographystyle{abbrv}
\bibliography{references}
\end{document}